\documentclass[11pt,a4paper]{article}
\usepackage{amsmath, amssymb}
\usepackage{graphicx,exscale,times,mathrsfs,a4wide,theorem,fancyheadings}
\numberwithin{equation}{section}
\newcommand{\wedg}{\mathbin{\scriptstyle{\wedge}}}
\newcommand{\lctr}{\mathbin{\lrcorner}}
\newcommand{\contr}{\lctr}
\newcommand{\Wedge}{\Lambda}
\newcommand{\Range}{\operatorname{im}}
\newcommand{\Nullspace}{\operatorname{ker}}
\renewcommand{\to}{\mathop{\,\rightarrow\,}\limits}
\newcommand{\supp}{\operatorname{supp}}

\let\curl=\rot
\newcommand {\grad}{\operatorname{\mathbf{grad}}}
\newcommand {\Div}{\operatorname{div}}

\newcommand{\Cinf}{{\mathscr C}^\infty}
\let\ov\overline

\newcommand  {\R}{\mathbb{R}}
\newcommand  {\N}{\mathbb{N}}

\let\te\textstyle

\let\cal\mathscr
\newtheorem{theorem}{Theorem}[section]
\newtheorem{lemma}[theorem]{Lemma}  
\newtheorem{proposition}[theorem]{Proposition}
\newtheorem{corollary}[theorem]{Corollary}
\newtheorem{remark}[theorem]{Remark}
\newenvironment{Remark}{\begin{remark}
      \parindent=0pt \parskip=0.5ex \rm }{\hfill $\triangle$\end{remark}}
\newtheorem{definition}[theorem]{Definition}

\newenvironment{proof}{\parindent=0pt \parskip=0.5ex 
   \abovedisplayskip = 0.8 \abovedisplayskip
   \belowdisplayskip=\abovedisplayskip
   \textbf{Proof: \ }}{\hfill $\square$\smallskip}

\title{On Bogovski\u{\i} and regularized Poincar\'e integral operators for 
 de Rham complexes on Lipschitz domains}

\author{Martin Costabel \and Alan McIntosh}
\date{December 29, 2008}
\begin{document}

\maketitle

\begin{abstract}
We study integral operators related to a regularized version of the
classical Poincar\'e path integral and the adjoint class generalizing Bogovski\u{\i}'s integral operator, acting on differential forms in $\R^n$. We prove that these operators are pseudodifferential operators of order $-1$. The Poincar\'e-type operators map polynomials to polynomials and can have applications in finite element analysis. 

For a domain starlike with respect to a ball, the special support properties of the operators imply regularity for the de Rham complex without boundary conditions (using Poincar\'e-type operators) and with full Dirichlet boundary conditions (using Bogovski\u{\i}-type operators). For bounded Lipschitz domains, the same regularity results hold, and in addition we show that the cohomology spaces can always be represented by $\Cinf$ functions.

\medskip

\noindent
\emph{2000 Mathematics Subject Classification.} 
Primary 35B65, 35C15; 
Secondary 58J10, 47G30 

\noindent
\emph{Key words and phrases.}
Exterior derivative, differential forms, Lipschitz domain, Sobolev spaces, pseudodifferential operator

\end{abstract} 

\section{Introduction}

In \cite{Bogovski79}, Bogovski\u{\i} introduced an integral operator $T$ with two remarkable properties: \\
- If $f$ is a function satisfying $\int f(x)dx=0$, then $u=Tf$ solves the partial differential equation $\Div u=f$, and\\
- If the bounded domain $\Omega\subset\R^n$ is starlike with respect to an open ball $B$, then $T$ maps the Sobolev space $W^{m-1,p}_0(\Omega)$ boundedly to $W^{m,p}_0(\Omega)^n$ for all $m\ge0$ and $1<p<\infty$.

This implies for a large class of domains $\Omega$, including all bounded Lipschitz domains, the solvability in $W^{m,p}_0(\Omega)^n$ of the equation 
$\Div u=f$ for $f\in W^{m-1,p}_0(\Omega)$ satisfying the integrability condition $\int f dx=0$. This means that there is no loss of regularity, and the support is preserved. 

This operator is now a classical tool in the theory of the equations of hydrodynamics \cite{Galdi94I}. It was recently noticed that its range of continuity can be extended to Sobolev spaces of negative order of regularity \cite{GeissertHeckHieber06}, and the study of more refined mapping properties has been instrumental in obtaining sharp regularity estimates for powers of the Stokes operator \cite{MitreaMonniaux_JFA08}.

Bogovski\u{\i}'s integral operator $T$ makes use of a smoothing function
\begin{equation}
\label{eq:theta}
 \theta\in \Cinf_0(\R^n)\,, \quad
 \supp\theta\subset B\,, \quad
 \int \theta(x)\,dx = 1
\end{equation}
when $\Omega$ is starlike with respect to an open ball $B$,
and is defined by
\begin{equation}
\label{eq:bogo}
   T f(x) = \int_\Omega f(y) \frac{x-y}{|x-y|^n}
              \int_{|x-y|}^\infty \theta\Big(y+ r\frac{x-y}{|x-y|}\Big) \,r^{n-1} \,dr \,dy\;.
\end{equation}
Applying the change of variables 
$(y,r)\mapsto(a,t)=(x+r\frac{y-x}{|x-y|},1-\frac{|x-y|}{r})$, one sees that the formally adjoint integral operator $T\,'$ is given by a smoothed-out path integral which defines the potential $v=T\,'u$ of a conservative vector field $u$, thus giving a solution of the equation $\grad v=u$:
\begin{equation}
\label{eq:bogoprime}
  T\,'u(x)=-\int \theta(a) J_a u(x)\,da\;, \quad
  J_a u(x) = (x-a)\cdot \int_0^1 u\bigl(a+t(x-a)\bigr)\,dt\;.
\end{equation}
The standard proof of Poincar\'e's lemma in differential geometry via ``Cartan's magic formula'' \cite[Theorem 13.2]{Taylor96I} uses a generalization of the path integral $J_a$ in \eqref{eq:bogoprime} to construct a right inverse of the exterior derivative operator for closed differential forms. A typical example in $\R^3$ is the path integral
\begin{equation}
\label{poincare3d}
 R_a u(x) = - (x-a)\times \int_0^1 u\bigl(a+t(x-a)\bigr)\, t \,dt
\end{equation} 
which provides a solution of the equation $\curl v=u$ for a divergence-free vector field $u$. Under the name ``Poincar\'e map'', this integral operator has recently been used in the analysis of finite element methods for Maxwell's equations \cite{GopalaDemko04,CoDaDe08}. Three properties of the operator $R_a$
are important for this application:\\
- $R_a$  maps polynomial vector fields to polynomial vector fields\\
- If $\Omega$ is starlike with respect to $a$, then the restriction of $R_au$ to $\Omega$ depends only on the restriction of $u$  to $\Omega$\\
- $R_a$ maps $L^2(\Omega)^3$ boundedly to itself.

One of the results of the present paper is that the regularized version $R$ of $R_a$, given by
$$
  R u(x) = \int \theta(a) R_a u(x)\,da\;,
$$
while still preserving polynomials and the local domain of influence, defines a bounded operator from $W^{s,p}(\Omega)$ to $W^{s+1,p}(\Omega)$ for all $s\in\R$ and $1<p<\infty$, if $\Omega$ is starlike with respect to the ball $B$. Such an operator was used in Section 4 of \cite{Axelsson} to obtain an inverse to the exterior derivative operator in $L^2$ spaces.

In \cite{Mitrea_Duke04}, Mitrea studied the generalization of both the Bogovski\u{\i}-type and the regularized Poincar\'e-type integral operators acting on differential forms with coefficients in Besov or Triebel-Lizorkin spaces. In \cite{MitreaMitreaMonniaux08}, Mitrea, Mitrea and Monniaux extended this analysis to show that these operators are regularizing of order one on a large class of such function spaces and to obtain sharp regularity estimates for the ``natural'' boundary value problems of the exterior derivative operator on Lipschitz domains. There the non-smoothness of the boundary of the domain implies that the solutions of these boundary value problems are singular, and therefore the solution operator is bounded for certain intervals of the regularity index $s$ depending on the exponent $p$, whereas for certain critical indices the boundary value problem does not define an operator with closed range.

In this paper, we prove that the Bogovski\u{\i}-type and the regularized Poincar\'e-type integral operators are classical pseudodifferential operators of order $-1$ with symbols in the H\"ormander class $S^{-1}_{1,0}(\R^n)$. As is well known \cite[Chapter 6]{Triebel92}, this implies immediately that the operators act as bounded operators in a wide range of function spaces including H\"older, Hardy or Sobolev spaces, or more generally the Besov spaces $B_{pq}^s$ for $0<p,q\le\infty$, and the Triebel-Lizorkin spaces $F_{pq}^s$ for $0<p<\infty$, $0<q\le\infty$. In each case, the operators map differential forms with coefficients of regularity $s$ boundedly to differential forms of regularity $s+1$ and, if $\Omega$ is bounded and starlike with respect to a ball, the Bogovski\u{\i}-type operators act between spaces of distributions with compact support in $\ov\Omega$, and the Poincar\'e-type operators act between spaces of restrictions to $\Omega$.

As a consequence, we obtain regularity results for the exterior derivative operator on bounded Lipschitz domains, either in spaces with compact support, or in spaces without boundary conditions, and these regularity results hold without restriction on the regularity index $s$. In particular, we show that the cohomology spaces of the de Rham complex on a bounded Lipschitz domain, either with compact support, or without boundary conditions, can be represented independently of the regularity index $s$ by finite dimensional spaces of differential forms with $\Cinf$ coefficients.

Thus, by the end of the paper, we will have employed the Bogovski\u{\i}-type and the regularized Poincar\'e-type integral operators to construct finite dimensional spaces   ${\cal H}_\ell(\ov\Omega) \subset \Cinf(\ov\Omega,\Wedge^\ell)$ and ${\cal H}_{\ov\Omega,\ell}(\R^n) \subset \Cinf_{\ov\Omega}(\R^n,\Wedge^\ell)$, each independent of the degree of regularity $s$,  such that all of the following direct sum decompositions hold true. To do this we use finitely many coverings of $\ov\Omega$, each by finitely many starlike domains. (A similar procedure would work for a Lipschitz domain in a compact Riemannian manifold.) See the next section for definitions.

\begin{theorem}\label{intro}  
Let $\Omega$ be a bounded Lipschitz domain in $\R^n$, and let $0\le\ell\le n$. Then for the spaces without boundary conditons,
\begin{align*}
 \Nullspace \Bigl(d:\Cinf(\ov\Omega,\Wedge^{\ell})\to \Cinf(\ov\Omega,\Wedge^{\ell+1})\Bigr)
  &=
  d\, \Cinf(\ov\Omega,\Wedge^{\ell-1}) \,\oplus\, {\cal H}_\ell(\ov\Omega)\;,\\
 \Nullspace \Bigl(d:H^{s}(\Omega,\Wedge^{\ell})\to H^{s-1}(\Omega,\Wedge^{\ell+1})\Bigr) &=
 d\, H^{s+1}(\Omega,\Wedge^{\ell-1}) \, \oplus \, {\cal H}_\ell(\ov\Omega)\,\end{align*}
where the $H^s$ ($-\infty<s<\infty$) denote Sobolev spaces, and, more generally,
\begin{align*}
  \Nullspace \Bigl(d:B^{s}_{pq}(\Omega,\Wedge^{\ell})\to B^{s-1}_{pq}(\Omega,\Wedge^{\ell+1})\Bigr) &=
 d\, B^{s+1}_{pq}(\Omega,\Wedge^{\ell-1}) \, \oplus \, {\cal H}_\ell(\ov\Omega)\;,\\
 \Nullspace \Bigl(d:F^{s}_{pq}(\Omega,\Wedge^{\ell})\to F^{s-1}_{pq}(\Omega,\Wedge^{\ell+1})\Bigr) &=
 d\, F^{s+1}_{pq}(\Omega,\Wedge^{\ell-1}) \, \oplus \, {\cal H}_\ell(\ov\Omega)
\end{align*}
where the $B^{s}_{pq}$ ($-\infty<s<\infty$, $0<p,q\le\infty$)
denote Besov spaces, and the \\
$F^{s}_{pq}$ ($-\infty<s<\infty$, $0<p<\infty$, $0<q\le\infty$) denote Triebel-Lizorkin spaces.

For the spaces with compact support, and the same values of $s,p$ and $q$, we have
\begin{align*}
\Nullspace \Bigl(d:\Cinf_{\ov\Omega}(\R^n,\Wedge^{\ell})\to \Cinf_{\ov\Omega}(\R^n,\Wedge^{\ell+1})\Bigr)
  &=
  d\, \Cinf_{\ov\Omega}(\R^n,\Wedge^{\ell-1}) \,\oplus\, {\cal H}_{\ov\Omega,\ell}(\R^n)\;,\\ \Nullspace \Bigl(d:H^{s}_{\ov\Omega}(\R^n,\Wedge^{\ell})\to H^{s-1}_{\ov\Omega}(\R^n,\Wedge^{\ell+1})\Bigr) 
  &= 
 d\, H^{s+1}_{\ov\Omega}(\R^n,\Wedge^{\ell-1}) \, \oplus \, {\cal H}_{\ov\Omega,\ell}(\R^n)\;,\\
\Nullspace \Bigl(d:B^{s}_{pq\,\ov\Omega}(\R^n,\Wedge^{\ell})\to B^{s-1}_{pq\,\ov\Omega}(\R^n,\Wedge^{\ell+1})\Bigr) 
  &= 
 d\, B^{s+1}_{pq\,\ov\Omega}(\R^n,\Wedge^{\ell-1}) \, \oplus \, {\cal H}_{\ov\Omega,\ell}(\R^n)\;,\\
\Nullspace \Bigl(d:F^{s}_{pq\,\ov\Omega}(\R^n,\Wedge^{\ell})\to F^{s-1}_{pq\,\ov\Omega}(\R^n,\Wedge^{\ell+1})\Bigr) 
  &= 
 d\, F^{s+1}_{pq\,\ov\Omega}(\R^n,\Wedge^{\ell-1}) \, \oplus \, {\cal H}_{\ov\Omega,\ell}(\R^n)
  \;. 
 \end{align*}
\end{theorem}

We remark without further discussion that this result has applications for the local Hardy spaces 
$h^1_r(\Omega,\Lambda^\ell)=F^0_{12}(\Omega,\Lambda^\ell)$ and 
$h^1_z(\Omega,\Lambda^\ell)= F^0_{12\,\ov\Omega}(\R^n,\Lambda^\ell)$.

\section{Notation and definitions} 
\label{S:notation}
For a bounded domain $\Omega$ in $\R^n$, we consider four spaces of infinitely differentiable functions. Besides $\Cinf(\Omega)$, the space of all infinitely differentiable functions in $\Omega$, and $\Cinf_0(\Omega)$, the functions with compact support in $\Omega$, we also use the space of restrictions to $\Omega$
$$
  \Cinf(\ov\Omega) = 
    \{u\in\Cinf(\Omega) \mid \exists \tilde{u} \in \Cinf(\R^n):\; 
          u=\tilde{u} \text{ on }\Omega \}
$$
and the space of functions with support in $\ov\Omega$
$$
  \Cinf_{\ov\Omega}(\R^n) =
    \{ u\in\Cinf(\R^n) \mid \supp u \subset\ov\Omega \}\;.
$$
Thus $\Cinf(\ov\Omega)$ is a quotient space of $\Cinf(\R^n)$ (or $\Cinf_0(\R^n)$) modulo functions vanishing on $\Omega$, and $\Cinf_{\ov\Omega}(\R^n)$ is a subspace of $\Cinf(\R^n)$ (or $\Cinf_0(\R^n)$). Likewise, for functions or distributions of regularity $s\in\R$, we consider spaces of restrictions to $\Omega$ and spaces with compact support in $\ov\Omega$.

By the term \emph{bounded Lipschitz domain} $\Omega$ in $\R^n$ we mean a connected bounded open set which is strongly Lipschitz in the sense that in the neighborhood of each point of $\ov\Omega=\Omega\cup\partial\Omega$ it is congruent to the domain below the graph of a scalar Lipschitz continuous function of $n-1$ variables. 

A domain $\Omega$ is \emph{starlike} with respect to a set $B$ if for every $x\in\Omega$ the convex hull of $\{x\}\cup B$ is contained in $\Omega$. From the definitions, it is not hard to see that a bounded domain which is starlike with respect to an open ball is Lipschitz, and that conversely, every bounded Lipschitz domain is the union of a finite number of domains, each of which is starlike with respect to an open ball. 
For the latter, one can choose, for example, domains congruent to the domain below the graph of a Lipschitz continuous function of Lipschitz constant $L$, bounded below by $H>0$, defined on a ball of radius $R$ in $\R^{n-1}$. Such a domain will be starlike with respect to an open ball centered at the origin as soon as $RL<H$. 

To keep the notation simple, we use the Sobolev space $H^s=W^{s,2}$ as representative for a space of regularity $s$. But, as already mentioned, many of the following arguments remain valid if the $L^2$-based Sobolev space $H^s$ is replaced by the Sobolev-Slobodeckii space $W^{s,p}$ or the Bessel potential space $H^s_p$ ($1<p<\infty$) or, more generally, by any of $B_{pq}^s$ ($0<p,q\le\infty$) or $F_{pq}^s$ ($0<p<\infty$, $0<q\le\infty$).

We let $H^s(\Omega)$ denote the quotient space of $H^s(\R^n)$ by the subspace of distributions vanishing in $\Omega$, while we let $H^s_{\ov\Omega}(\R^n)$ denote the subspace of $H^s(\R^n)$ consisting of all distributions with support in $\ov\Omega$. 
Thus $H^s(\Omega)$, for which also equivalent intrinsic definitions exist, can be considered as a space of distributions on $\Omega$, whereas $H^s_{\ov\Omega}(\R^n)$ is a space of distributions on $\R^n$. 

Let us mention some well-known properties of these spaces that hold if $\Omega$ is a bounded Lipschitz domain. Proofs (for the spaces $W^{s,p}$, $s\in\R$, $1<p<\infty$) can be found in \cite[Chapter 1]{Grisvard85}: 
The intersection of all $H^s(\Omega)$, $s\in\R$, is $\Cinf(\ov\Omega)$ and the union of all $H^s(\Omega)$ is the space of all distributions on $\Omega$ that allow an extension to a neighborhood of $\ov\Omega$. 
Likewise, the intersection of all $H^s_{\ov\Omega}(\R^n)$ is $\Cinf_{\ov\Omega}(\R^n)$ and the union of all $H^s_{\ov\Omega}(\R^n)$ is the space of all distributions on $\R^n$ with support in $\ov\Omega$.
It is also well known that $H^s_{\ov\Omega}(\R^n)$, for which also Triebel's notation 
$\widetilde H^s(\Omega)$ is commonly used, can be identified with the space $H^s_0(\Omega)$, the closure of $\Cinf_0(\Omega)$ in $H^s(\Omega)$, if $s$ is positive and $s-\frac12$ is not an integer. For any $s\in\R$, $H^s_{\ov\Omega}(\R^n)$ is the closure of $\Cinf_0(\Omega)$ in $H^s(\R^n)$.
In our Hilbert space setting, for all $s\in\R$ the space $H^s_{\ov\Omega}(\R^n)$ is in a natural way isomorphic to the dual space of 
$H^{-s}(\Omega)$.

For differential forms we use standard notation which is, for example, defined in \cite{Sternberg83,Taylor96I}. The exterior algebra of $\R^n$ is $\Wedge^\ell$,
$0\le\ell\le n$,  where $\Wedge^0$ and $\Wedge^1$ are identified with
$\R$ and $\R^n$, respectively, and we set $\Wedge^\ell=\{0\}$ if $\ell<0$ or $\ell>n$.

Differential forms of order $\ell$ with coefficients in $H^s$ are denoted by $H^s(\Omega,\Wedge^\ell)$ and $H^s_{\ov\Omega}(\R^n,\Wedge^\ell)$. With the exterior derivative $d$ satisfying $d\circ d=0$ we then have the 
\emph{de Rham complex without boundary conditions}
\begin{equation}
\label{eq:dRwobc}
 0 \to H^s(\Omega,\Wedge^0) \to^d H^{s-1}(\Omega,\Wedge^1) \to^d\cdots
 \to^d H^{s-n}(\Omega,\Wedge^n) \to 0
\end{equation}
and the
\emph{de Rham complex with compact support}
\begin{equation}
\label{eq:dRwcs}
 0 \to H^s_{\ov\Omega}(\R^n,\Wedge^0) \to^d H^{s-1}_{\ov\Omega}(\R^n,\Wedge^1) \to^d\cdots
 \to^d H^{s-n}_{\ov\Omega}(\R^n,\Wedge^n) \to 0
\end{equation}

Besides these complexes we also consider the \emph{extended} de Rham complexes
without boundary conditions
\begin{equation}
\label{eq:edRwobc}
 0 \to \R \to^\iota H^s(\Omega,\Wedge^0) \to^d H^{s-1}(\Omega,\Wedge^1) \to^d\cdots
 \to^d H^{s-n}(\Omega,\Wedge^n) \to 0
\end{equation}
and 
with compact support
\begin{equation}
\label{eq:edRwcs}
 0 \to H^s_{\ov\Omega}(\R^n,\Wedge^0) \to^d H^{s-1}_{\ov\Omega}(\R^n,\Wedge^1) \to^d\cdots
 \to^d H^{s-n}_{\ov\Omega}(\R^n,\Wedge^n) \to^{\iota^*} \R \to 0
\end{equation}

Here the mapping denoted by $\iota$ in \eqref{eq:edRwobc} is the natural inclusion of constant functions, and 
$\iota^*$ in \eqref{eq:edRwcs} is the generalization to distributional coefficients with compact support of the integral
$u\mapsto \iota^*u=\int_{\R^n} u$ for an $n$-form $u$ with integrable coefficients.

The extended de Rham complexes \eqref{eq:edRwobc} and \eqref{eq:edRwcs} are exact at the left end because $\Omega$ is connected, and their exactness at the right end is the subject of Bogovski\u{\i}'s theorem mentioned in the introduction.
We will show in Section~\ref{S:RegdeRham} below that for bounded domains starlike with respect to a ball, both complexes \eqref{eq:edRwobc} and \eqref{eq:edRwcs} are exact for any $s\in\R$, and that for bounded Lip\-schitz domains both complexes \eqref{eq:dRwobc} and  \eqref{eq:dRwcs} have finite dimensional cohomology spaces whose dimension does not depend on $s$.

We will make use of the following standard algebraic operations in the exterior algebra which then also extend as pointwise operations to differential forms on domains of $\R^n$:

\begin{tabular}{lrcl}
the exterior product:& 
  $\wedg$&$:$& $\Wedge^\ell\times\Wedge^m \to \Wedge^{\ell+m}$\\
the interior product or contraction:&
  $\contr$&$:$& $\Wedge^\ell\times\Wedge^m \to \Wedge^{m-\ell}$ \\
the euclidean inner product:&
  $\langle a,b\rangle$&$:$&$\Wedge^\ell\times\Wedge^\ell \to \R$\\
the Hodge star operator:&
  $\star$&$:$& $\Wedge^\ell \to \Wedge^{n-\ell}$
\end{tabular}

We now give a list of well-known properties of these operations which will be sufficient for verifying the arguments used in our proofs below.

In particular we need the exterior product and the contraction with a vector $a\in\R^n$, identified with a $1$-form. For $a=(a_1,\dotsc,a_n)$ and 
$u=dx_{j_1}\wedg\dots\wedg dx_{j_\ell}$ with $j_1<\dots<j_\ell$, the contraction is given by 
$$
  a\contr u 
   = \sum_{k=1}^\ell 
     (-1)^{k-1} a_{j_k} dx_{j_1}\wedg\dots\wedg\widehat{dx}_{j_k}\wedg\dots\wedg dx_{j_\ell}
$$
where the notation $\widehat{dx}_{j_k}$ means that the corresponding factor is to be omitted.
In the special case of $\R^3$, this corresponds to the following classical operations of vector algebra:
\begin{align*}
  u \text{ scalar, interpreted as $0$-form:} \qquad a\wedg u &= ua & a\contr u &=0\\  
  u \text{ scalar, interpreted as $3$-form:} \qquad a\wedg u &= 0  & a\contr u &=ua\\  
  u \text{ vector, interpreted as $1$-form:} 
     \qquad a\wedg u &= a\times u & a\contr u &=a\cdot u\\  
  u \text{ vector, interpreted as $2$-form:} 
     \qquad a\wedg u &= a\cdot u & a\contr u &= - a\times u  
\end{align*}
Some useful formulas for $u,v\in\Wedge^\ell$, $w\in\Wedge^{\ell+1}$, $a\in\Wedge^1$ are:
\begin{align}
\label{eq:starstar}
  \star \star u &=  (-1)^{\ell(n-\ell)}u\\
\label{eq:starwedge}
  \star (a\wedg u) &= (-1)^{\ell}\, a\contr (\star u)\\
\label{eq:starprod}
  \langle u,v\rangle &= \star (u\wedg\star v) = \langle \star u,\star v\rangle\\
\label{eq:prodcontr}
  \langle w,a\wedg u\rangle &= \langle u,a\contr w\rangle\\
\intertext{We note the product rule of the exterior derivative for an $\ell$-form $u$ and an $m$-form $v$}
\label{eq:prodrule}
 d(u\wedg v) &= (du)\wedg v + (-1)^\ell u\wedg(dv) \,.
\end{align}

Finally, with the $L^2$ scalar product for $\ell$-forms $u$ and $v$,
$$
  (u,v) = \int_\Omega \langle u(x), v(x)\rangle\, dx
$$
and the co-derivative $\delta$, there holds
\begin{gather}
\label{eq:deltad}
  (\delta u,v) = (u,dv)\;,\\
\label{eq:dstardelta}
  \star \,\delta \,=\, (-1)^{\ell}\, d\star  \quad\mbox{ and }\quad
  \star d \,=\, (-1)^{\ell-1} \,\delta\star  \quad\mbox{ on $\ell$-forms}\;.
\end{gather}

\section{The Bogovski\u{\i} and Poincar\'e integral operators}
\label{S:BogoPoin}

In this section, we fix a function 
$\theta\in\Cinf_0(\R^n)$ with support in a ball $B$ satisfying $\int\!\theta(x)\,dx=1$.

\subsection{Definition, support properties}

For 
$\ell\in\{0,\dotsc,n\}$, define the kernel $G_\ell$ by
\begin{equation}
\label{eq:Gl}
  G_\ell(x,y) = 
  \int_1^\infty (t-1)^{n-\ell}t^{\ell-1}\theta \bigl(y+t(x-y)\bigr)\,dt \;.
\end{equation}
\begin{definition}
\label{D:4iops}
 For a differential form $u\in\Cinf_0(\R^n,\Wedge^\ell)$, define two integral operators:
\begin{align}
\label{eq:Rl}
  R_\ell u(x) &=
  \int G_{n-\ell+1}(y,x) \,(x-y)\contr u(y)\,dy
  \quad (1\le\ell\le n)\\
\label{eq:Tl}
  T_\ell u(x) &= 
  \int G_{\ell}(x,y) \,(x-y)\contr u(y)\,dy
  \quad (1\le\ell\le n)
\end{align}
We  refer to $R_\ell$ as Poincar\'e-type operators, and to $T_\ell$ as Bogovski\u{\i}-type operators.
\end{definition}

In order to see that the integrals in Definition~\ref{D:4iops} exist, we rewrite the kernel $G_\ell$:
\begin{align}
\label{eq:Glsumhomog}
  G_\ell(x,y) &= 
     \int_0^\infty \tau^{n-\ell} (\tau+1)^{\ell-1}\theta\bigl(x+\tau(x-y)\bigr)\,d\tau\notag\\
   &= \sum_{k=0}^{\ell-1}\tbinom{\ell-1}{k} 
     \int_0^\infty \tau^{n-k-1} \theta\bigl(x+\tau(x-y)\bigr)\,d\tau\notag\\
   &= \sum_{k=0}^{\ell-1}\tbinom{\ell-1}{k} \,
     |x-y|^{k-n} \int_0^\infty r^{n-k-1} \theta\bigr(x+r\frac{x-y}{|x-y|}\bigr)\,dr \;.
\end{align}
This representation as a finite sum of homogeneous functions gives a bound 
\begin{equation}
\label{eq:Glweaksing}
  |G_\ell(x,y)\,(x-y)|\le C(x)\, |x-y|^{-n+1}\,,
\end{equation}
where $C(x)$ depends on $\|\theta\|_{L^\infty}$ and the size of the ball $B$, and is uniformly bounded for $x$ in a bounded set. 
Hence the integrals in Definition~\ref{D:4iops} are weakly singular and therefore convergent.

As one can readily see from the definitions, the two integral operators are related by duality: If we introduce operators $Q_\ell$ and $S_\ell$ by Hodge star duality, so that for $0\le\ell\le n-1$ and $u\in\Cinf_0(\R^n,\Wedge^\ell)$
\begin{equation}
\label{eq:starQRST}
  \star Q_\ell\, u = (-1)^{\ell-1}\,R_{n-\ell}\,(\star u) \quad \text{ and }\quad  
  \star S_\ell\, u = (-1)^{\ell-1}\,T_{n-\ell}\,(\star u)\,,
\end{equation}
then we have for $v\in\Cinf_0(\R^n,\Wedge^{\ell+1})$
\begin{equation}
\label{eq:adjQRST}
   \bigl(v,Q_\ell u\bigr) = \bigl(T_{\ell+1}v, u\bigr) \quad \text{ and }\quad
   \bigl(v,S_\ell u\bigr) = \bigl(R_{\ell+1}v, u\bigr) \;.
\end{equation}
Denoting the formal adjoint operator with respect to the $L^2$ duality by a prime, we have therefore
\begin{equation}
\label{eq:dualRT}
 \star R_{\ell} = (-1)^\ell\, T\,'_{\!n-\ell+1}\,\star
\end{equation}
In order to see other properties of the operators, we apply a different change of variables. Let us write this in detail for the operator $R_\ell$. We use the change of variables
$a=x+t(y-x)$ and then replace $(t-1)/t$ by $t$.
\begin{align}
\label{eq:regPoincare} 
 R_\ell u(x) &= 
   \int \int_1^\infty (t-1)^{\ell-1}t^{n-\ell}
          \theta\bigl(x+t(y-x)\bigr)\, (x-y)\contr u(y) \,dt\,dy\notag\\
  &=
   \int \int_1^\infty (t-1)^{\ell-1}t^{-\ell-1}\theta(a) \,(x-a)\contr
      u\bigl(x+(a-x)/t\bigr)\,dt\,da\notag\\
  &=
   \int \theta(a) \,(x-a)\contr \int_0^1 t^{\ell-1}\,
      u\bigl(a+t(x-a)\bigr)\,dt\,da\;.
\end{align}
From this form of $R_\ell$, one sees immediately that it maps differential forms with polynomial coefficients to differential forms with polynomial coefficients and also
$\Cinf(\R^n,\Wedge^\ell)$ to $\Cinf(\R^n,\Wedge^{\ell-1})$,
and that $R_\ell u(x)$ depends only on the values of $u$ in the convex hull of $B\cup\{x\}$, that is, the starlike hull of $\{x\}$ with respect to the ball $B$. This implies in particular that if $\Omega$ is open and starlike with respect to $B$, then $R_\ell$ maps 
$\Cinf(\Omega,\Wedge^\ell)$ to $\Cinf(\Omega,\Wedge^{\ell-1})$ and also
$\Cinf(\ov\Omega,\Wedge^\ell)$ to $\Cinf(\ov\Omega,\Wedge^{\ell-1})$.

Rewriting $T_\ell$ in the same way, we get
\begin{equation}
\label{eq:Bogo}
 T_\ell u(x) = -
   \int \theta(a) \,(x-a)\contr \int_1^\infty t^{\ell-1}\,
      u\bigl(a+t(x-a)\bigr)\,dt\,da\;.
\end{equation}
From this form of $T_\ell$, because of the unbounded interval of integration in $t$, one cannot immediately conclude that $T_\ell$ maps $\Cinf$ functions to $\Cinf$ functions.
But if $u\in\Cinf_0(\R^n,\Wedge^\ell)$, one sees that $T_\ell u$ is $\Cinf$ on 
$\R^n\setminus\supp\theta$, and that 
$T_\ell u(x)=0$ unless $x$ lies in the starlike hull of $\supp u$ with respect to $B$. 
Thus if $\Omega$ is open and starlike with respect to $B$, then   
$u\in\Cinf_0(\Omega,\Wedge^\ell)$ implies $\supp T_\ell u \subset\Omega$, and, if $\Omega$ is bounded, then
$u\in\Cinf_{\ov\Omega}(\R^n,\Wedge^\ell)$ implies $\supp T_\ell u \subset\ov\Omega$. 
The fact that $T_\ell$ indeed maps 
$\Cinf_0(\R^n,\Wedge^\ell)$ to $\Cinf_0(\R^n,\Wedge^{\ell-1})$ will be a consequence of Theorem~\ref{T:pseudo} below.

\subsection{Homotopy relations}

Cartan's formula for the Lie derivative of a differential form with respect to a vector field can be written as
$$
  \frac{d}{dt}F^*_t u = F^*_t \bigl(d(X_t\contr u) + X_t\contr du\bigr)\;,
$$
where $F^*_t$ denotes the pull-back by the flow $F_t$ associated with the vector field $X_t$. Here we consider the special case of the dilation flow with center $a$
$$
  F_t(x) = a+t(x-a) \quad\text{ with vector field } X_t=x-a\;,
$$
which gives a pull-back of
$$
  F^*_t u(x) = t^\ell \,u\bigl(a+t(x-a)\bigr) \quad \text{for an $\ell$-form } u \;.
$$
This leads to the formula
\begin{equation}
\label{eq:Cartan}
 \frac{d}{dt}(t^\ell u\bigl(a+t(x-a)\bigr) =
    d\Bigl(t^{\ell-1}(x-a)\contr u\bigl(a+t(x-a)\bigr)\Bigr) +
    t^{\ell}(x-a)\contr du\bigl(a+t(x-a)\bigr)
\end{equation}
which can also be verified elementarily from the formulas we gave in Section~\ref{S:notation}.

Integrating \eqref{eq:Cartan} from $0$ to $1$ and comparing with \eqref{eq:regPoincare}, 
we find the homotopy relations, valid for all $u\in\Cinf_0(\R^n,\Wedge^\ell)$ 
\begin{equation}
\begin{aligned}
\label{eq:dR+Rd=1}
 dR_\ell u + R_{\ell+1} du &= u\; &&(1\le\ell\le n-1)\;;\\
 R_1 du &= u - \bigl(\theta,u\bigr)\; &&(\ell=0)\;;\\
 dR_nu &=u &&(\ell=n)\;.
\end{aligned}
\end{equation}
One could be tempted to integrate Cartan's formula from $1$ to $\infty$ and compare with \eqref{eq:Bogo}, thus formally obtaining a similar homotopy relation for $T_\ell$ directly. The result is indeed true except for $\ell=n$, but for a rigorous proof we prefer to use the duality relation \eqref{eq:dualRT} to deduce corresponding anticommutation relations for $T_\ell$ from the relations \eqref{eq:dR+Rd=1} which are already proved. Here is what one obtains for $u\in\Cinf_0(\R^n,\Wedge^\ell)$:
\begin{equation}
\begin{aligned}
\label{eq:dT+Td=1}
  dT_\ell u + T_{\ell+1} du &= u\; &&(1\le\ell\le n-1)\;;\\
  T_1 du &= u \; &&(\ell=0)\;;\\
 dT_n u &=u - (\te\int\! u)\star\theta&&(\ell=n)\;.
\end{aligned}
\end{equation}
Here we consider $\theta$ as an element of $\Cinf_0(\R^n,\Wedge^0)$, so that for another $0$-form $u$ we have the $L^2$ scalar product $\bigl(\theta,u\bigr)=\int\theta(a)u(a)da$, and $\star\theta$ is the $n$-form $\theta(x)dx_1\wedg\dots\wedg dx_n$. 

The formulas for the endpoints $\ell=0$ and $\ell=n$ correspond to the two extended de Rham complexes without boundary conditions and with compact support, see \eqref{eq:edRwobc} and \eqref{eq:edRwcs}. To see this, let us extend the definition of the exterior derivative by writing $\ov d$ for all the mappings of the complex
$$
 0 \to\R\to^\iota\Cinf(\ov\Omega,\Wedge^0) \to^d \Cinf(\ov\Omega,\Wedge^1)
 \to^d\cdots \to^d \Cinf(\ov\Omega,\Wedge^n) \to 0
$$ 
and $\underline d$ for all the mappings of the complex
$$
 0 \to \Cinf_{\ov\Omega}(\R^n,\Wedge^0) \to^d \Cinf_{\ov\Omega}(\R^n,\Wedge^1) \to^d\cdots
 \to^d \Cinf_{\ov\Omega}(\R^n,\Wedge^n) \to^{\iota^*} \R \to 0
$$
where $\iota$ is the inclusion mapping for constant functions and 
$\iota^*=(\star\iota)'$ denotes the integral $u\mapsto\int\!u$ for $n$-forms.

If we correspondingly extend the definitions of $R_\ell$ and $T_\ell$ by
$$
\begin{aligned}
  R_0u &:=\bigl(\theta,u\bigr)\,\text{ for $0$-forms }u\,, &
  R_{n+1} &:=0\,,\\
  T_{n+1}u &:= \star(u\theta) \,\text{ for }u\in\R\,,  &
  T_0 &:= 0\,,
\end{aligned}
$$
then we can write the relations \eqref{eq:dR+Rd=1} and \eqref{eq:dT+Td=1} simply as
\begin{equation}
\label{eq:R&T&d}
 \ov d\,R_\ell u + R_{\ell+1} \,\ov d u  = u\quad \text{ and }\quad
 \underline d\,T_\ell u + T_{\ell+1} \,\underline du = u\quad \text{ for all }\;
  0\le\ell\le n.
\end{equation}

\subsection{Continuity}

The most important result about analytic properties of our integral operators is the following.
\begin{theorem}
\label{T:pseudo}
 The operators $R_\ell$ and $T_\ell$ defined in Definition~\ref{D:4iops} are pseudodifferential operators on $\R^n$ of order $-1$ with symbols in the H\"ormander symbol class $S^{-1}_{1,0}(\R^n)$.
\end{theorem}
\begin{proof}
 For basic facts about pseudodifferential operators, see for example \cite{Taylor81,Taylor96II,Wloka95}. We are using here the local symbol class 
$S^{-1}_{1,0}(\R^n)$ that consists of functions $a\in\Cinf(\R^n\times\R^n)$ satisifying for any compact set $M \subset\R^n$ and any multi-indices $\alpha,\beta\in \N_0^n$, estimates of the form
\begin{equation}
\label{eq:S^-1_1,0}
     |\partial_x^\alpha \partial_\xi^\beta  a(x,\xi)| \le
     C_{\alpha\beta}(M)\,(1+|\xi|)^{-1-|\beta|} \quad \forall (x,\xi)\in M\times\R^n\;.
\end{equation}
The proof will show that the constants $C_{\alpha\beta}$ are polynomially bounded in $x\in\R^n$, but this is not important here, since we are only interested in the local behavior. 

   We give the proof for the operator $T_\ell$. For $R_\ell$ the result then follows from \eqref{eq:dualRT} by applying the Hodge star operator which is a purely algebraic operation on basis vectors in the exterior algebra and does not change coefficients of differential forms, and by taking $L^2$ adjoints, which according to the calculus of pseudodifferential operators does not lead out of this class.

Thus we consider the integral operator defined by
$$
T_\ell u(x) = 
  \int G_{\ell}(x,y) \,(x-y)\contr u(y)\,dy
$$
with the kernel $G_\ell$ given in \eqref{eq:Gl}. Writing the differential forms in components, we see that for $j,\ell\in\{1,\dotsc,n\}$ we need to study the following operator $K$ acting on scalar functions $u$:
\begin{multline}
 Ku(x) = \int_{\R^n}k(x,x-y)\,u(y)\,dy\\
 \text{ with }\;
 k(x,z) = z_j \int_0^\infty s^{n-\ell}(s+1)^{\ell-1}\theta(x+sz)\,ds\;
 \text{ for } x,z\in\R^n\,.
\end{multline}
We write $k(x,z)=k_0(x,z)+k_1(x,z)$ with 
\begin{align*}
  k_0(x,z)&= z_j\, \int_0^1 s^{n-\ell}(s+1)^{\ell-1} \theta(x+sz)\,ds\,,\\
  k_1(x,z)&= z_j\, \int_1^\infty s^{n-\ell}(s+1)^{\ell-1} \theta(x+sz)\,ds\,.
\end{align*}
It is clear that $k_0\in\Cinf(\R^{2n})$, and therefore only $k_1$ needs to be analyzed.
If
$\supp\theta\subset B_\epsilon(0)$, then 
$$
  k_1(x,z) = 0 \quad\mbox{ for } |z|\ge|x|+\epsilon\,,
$$
and
we have already seen in \eqref{eq:Glweaksing} that $z\mapsto k_1(x,z)$ is weakly singular. It is therefore integrable over $\R^n$, so we can write
its Fourier transform as the convergent integral
\begin{align*}
  \hat k_1(x,\xi) &= 
    \int_{\R^n} e^{-i\langle\xi,z\rangle}k_1(x,z)\,dz\\
  &= \int_1^\infty s^{n-\ell}(s+1)^{\ell-1} 
     \int e^{-i\langle\xi,z\rangle} z_j \theta(x+sz)\,dz\,ds\,,
\end{align*}
and we can represent the operator $K$ as
\begin{equation*}
  Ku(x) \,=\,
   \int_{\R^n} k_0(x,x-y)\,u(y)\,dy +
   (2\pi)^{-n} \int_{\R^n} e^{i\langle\xi,x\rangle} 
   \hat k_1(x,\xi)\,\hat u(\xi)\,d\xi\;.
\end{equation*}

The proof will be complete once we show that the symbol $\hat k_1$ of the operator $K$ satisfies the estimates \eqref{eq:S^-1_1,0}, namely 
for any multi-indices $\alpha,\beta\in \N_0^n$ and 
$x,\xi\in\R^n$:
\begin{equation}
\label{eq:dadb}
  |\partial_x^\alpha \partial_\xi^\beta \hat k_1(x,\xi)| \le
   C_{\alpha\beta}(x)\,(1+|\xi|)^{-1-|\beta|}
\end{equation}
where $C_{\alpha\beta}(x)$ is bounded for $x$ in any compact set. 

With the change of variables $(t,y)=(1/s,x+sz)$ we can write
\begin{align}
\nonumber
  \hat k_1(x,\xi)  &= 
    \int_0^1(t+1)^{\ell-1}e^{it\langle\xi,x\rangle}
    \int e^{-it\langle\xi,y\rangle}\,(y_j-x_j)
    \theta(y) \,dy\,dt\\
\label{eq:k1hat}
  &=  \int_0^1(t+1)^{\ell-1}e^{it\langle\xi,x\rangle}
    \Bigl(i(\partial_j\hat\theta)(t\xi)-x_j\hat\theta(t\xi)\Bigr)\,dt\,.
\end{align}
Here $\hat\theta$ is the Fourier transform of $\theta\in\Cinf_0(\R^n)$, thus a rapidly decreasing $\Cinf$ function.
The representation \eqref{eq:k1hat}
shows that 
\begin{equation}
\label{eq:h1hatCinf}
  \hat k_1 \in \Cinf(\R^{2n}) \quad
  \text{ and } \;
  |\hat k_1(x,\xi)| \le C_\theta(1+|x|)
\end{equation}
where $C_\theta$ depends only on $\theta$. 
Writing $\tau=t|\xi|$ and $\omega=\xi/|\xi|$, we find 
\begin{equation}
\label{eq:h1hathom}
  \hat k_1(x,\xi) \,=\, |\xi|^{-1}
  \int_0^{|\xi|}\bigl(1+\frac\tau{|\xi|}\bigr)^{\ell-1}e^{i\tau\langle\omega,x\rangle}
    \Bigl(i(\partial_j\hat\theta)(\tau\omega)-x_j\hat\theta(\tau\omega)\Bigr)
   \,d\tau 
\end{equation}
and hence
$$
 |\hat k_1(x,\xi)| \le 
 |\xi|^{-1} 2^{\ell-1} \int_0^\infty\bigl(|\partial_j\hat\theta)(\tau\omega)|
   +|x_j\hat\theta(\tau\omega)|\bigr)\,d\tau \le
   (1+|x|)\,C_\theta\,|\xi|^{-1}\;.
$$
Thus we have shown \eqref{eq:dadb} for $|\alpha|=|\beta|=0$.

Similarly, by taking derivatives in \eqref{eq:k1hat}, we can write for any
multi-indices $\alpha,\beta$:

\begin{equation}
\label{eq:dadbk1hat}
  \partial_x^\alpha \partial_\xi^\beta \hat k_1(x,\xi) \,=\, 
  \int_0^1(t+1)^{\ell-1}e^{it\langle\xi,x\rangle}
    t^{|\beta|}
    \Bigl(p_{\alpha\beta}(x,t\xi,\partial)\hat\theta\Bigr)(t\xi)
   \,dt\;,
\end{equation}
where $p_{\alpha\beta}(x,\xi,\partial)$ is a partial differential
operator of order $|\beta|+1$ with polynomial coefficients of degree
$\le|\beta|+1$ in $x$ and $\le|\alpha|$ in $\xi$. We obtain an
immediate estimate
\begin{equation}
\label{eq:dadbat0}
  |\partial_x^\alpha \partial_\xi^\beta \hat k_1(x,\xi)| \le
   C_{\alpha\beta}(x)\,(1+|\xi|)^{|\alpha|}\;,
\end{equation}
and after the change of variables $\tau=t|\xi|$ with
$\omega=\xi/|\xi|$:
\begin{equation}
\label{eq:dadbk1hat2}
  \partial_x^\alpha \partial_\xi^\beta \hat k_1(x,\xi) \,=\, 
  \int_0^{|\xi|}\bigl(1+\frac\tau{|\xi|}\bigr)^{\ell-1}e^{i\tau\langle\omega,x\rangle}
    \tau^{|\beta|}
    \Bigl(p_{\alpha\beta}(x,\tau\omega,\partial)\hat\theta\Bigr)(\tau\omega)
   \,d\tau\,|\xi|^{-1-|\beta|}\;.
\end{equation}
This gives a second estimate
\begin{equation}
\label{eq:dadbatinf}
  |\partial_x^\alpha \partial_\xi^\beta \hat k_1(x,\xi)| \le
   C_{\alpha\beta}(x)\,|\xi|^{-1-|\beta|}\;.
\end{equation}
In \eqref{eq:dadbat0} and \eqref{eq:dadbatinf}, $C_{\alpha\beta}(x)$ is bounded for $x$ in any compact set. 
(One can see that
$
  C_{\alpha\beta}(x) \le C_{\alpha\beta,\theta}\cdot(1+|x|)^{1+|\beta|}
$
where $C_{\alpha\beta,\theta}$ depends only on $\alpha$, $\beta$ and
$\theta$.)

This shows \eqref{eq:dadb} and completes the proof.
\end{proof}

An immediate consequence of the theorem is that the two integral operators map differential forms with $\Cinf_0$ coefficients to differential forms with $\Cinf$ coefficients. Taking into account the support properties deduced above from the representations \eqref{eq:regPoincare} and \eqref{eq:Bogo}, we get the following statements, where we use the standard topologies for the function spaces. These statements follow also from the results in \cite[Theorem 4.1]{MitreaMitreaMonniaux08}.
\begin{corollary}
\label{C:QRSTCinf}
 The integral operators defined in Definition~\ref{D:4iops} define continuous mappings
\begin{align*}
 R_\ell&: \Cinf(\R^n,\Wedge^\ell)\to\Cinf(\R^n,\Wedge^{\ell-1})\,, 
 & T_\ell&: \Cinf_0(\R^n,\Wedge^\ell)\to\Cinf_0(\R^n,\Wedge^{\ell-1})\,.
\end{align*}
If $\Omega\subset\R^n$ is a bounded domain starlike with respect to a ball $B$ containing  $\supp\theta$, then the operators define continuous mappings
\begin{align*}
 R_\ell&: \Cinf(\Omega,\Wedge^\ell)\to\Cinf(\Omega,\Wedge^{\ell-1})\,, 
 & R_\ell&: \Cinf(\ov\Omega,\Wedge^\ell)\to\Cinf(\ov\Omega,\Wedge^{\ell-1})\,,\\
 T_\ell&: \Cinf_0(\Omega,\Wedge^\ell)\to\Cinf_0(\Omega,\Wedge^{\ell-1})\,, 
 & T_\ell&: \Cinf_{\ov\Omega}(\R^n,\Wedge^\ell)\to\Cinf_{\ov\Omega}(\R^n,\Wedge^{\ell-1})\,.
\end{align*}
\end{corollary}

Either by duality or by extension using standard continuity properties of pseudodifferential operators, the two operators can be defined on differential forms with distributional coefficients, in the case of the Poincar\'e-type operators $R_\ell$ for arbitrary distributions from ${\cal D}'(\R^n,\Wedge^\ell)$ and in the case of the Bogovski\u{\i}-type operators $T_\ell$ for distributions with compact support in $\R^n$.

For finite regularity, the standard continuity properties of pseudodifferential operators together with the support properties immediately imply results of the following type.

\begin{corollary}
\label{C:QRSTH^s}
Let $\Omega\subset\R^n$ be a bounded domain starlike with respect to a ball $B$ containing  $\supp\theta$. Then the two integral operators define bounded operators for any $s\in\R${\rm:}
\begin{align*}
 R_\ell&: H^s(\Omega,\Wedge^\ell)\to H^{s+1}(\Omega,\Wedge^{\ell-1})\,, 
 & T_\ell&: H^s_{\ov\Omega}(\R^n,\Wedge^\ell)\to H^{s+1}_{\ov\Omega}(\R^n,\Wedge^{\ell-1})\,.\end{align*}
\end{corollary}

\begin{Remark}
\label{R:QRSTH^s}
Corollary \ref{C:QRSTH^s} remains valid when $H^s$ is replaced by 
$B^s_{pq}\, (0<p\leq \infty, 0<q\leq \infty)$, or by 
$F^s_{pq}\, (0<p<\infty, 0<q\leq \infty)$.   
The spaces $B^s_{pq}(\Omega,\Wedge^\ell)$ and $F^s_{pq}(\Omega,\Wedge^\ell)$ are defined as quotient spaces, and the  spaces 
$B^s_{pq\,\overline\Omega}(R^n,\Wedge^\ell)$ and 
$F^s_{pq\,\overline\Omega}(R^n,\Wedge^\ell)$ are defined as subspaces, in an analogous way to the spaces 
$H^s(\Omega,\Wedge^\ell)$ and $H^s_{\overline\Omega}(R^n,\Wedge^\ell)$.  
They include the special cases of Sobolev spaces $W^{s,p}=F^s_{p,2}$, and local Hardy spaces $h^1_r(\Omega,\Wedge^\ell) = F^0_{1,2}(\Omega,\Wedge^\ell)$ and 
$h^1_z(\Omega,\Lambda^l) = F^0_{1,2\,\overline\Omega}(R^n,\Wedge^\ell)$. See Chapter 6 of \cite{Taylor96II}.
\end{Remark}

In all these cases, the commutation relations \eqref{eq:dR+Rd=1}--\eqref{eq:R&T&d} remain valid. What this implies for the regularity of the de Rham complex and its cohomology is the subject of the next section.

\section{Regularity of the de Rham complex}
\label{S:RegdeRham}

\subsection{Starlike domains}

The homotopy relations \eqref{eq:R&T&d} together with the mapping properties from Corollary~\ref{C:QRSTH^s} imply the existence of regular solutions of the equation $du=0$, as we now state. There are similar results in the $\Cinf$ spaces which follow from Corollary~\ref{C:QRSTCinf}.

\begin{proposition}
\label{P:starlike}
Let $\Omega \subset\R^n$ be a bounded domain, starlike with respect to a ball $B$. \\
{\rm(i)} For any $s\in\R$ and $\ell\in\{1,\dotsc,n\}$, let $u\in H^s(\Omega,\Wedge^\ell)$ satisfy $du=0$ in $\Omega$. Then there exists $v\in H^{s+1}(\Omega,\Wedge^{\ell-1})$ such that $dv=u$, and there is a constant $C$ independent of $u$ such that
$$
   \Vert v\Vert_{H^{s+1}(\Omega)} \le C\, \Vert u\Vert_{H^{s}(\Omega)}\;.
$$
For $\ell=n$ the condition $du=0$ is always satisfied.\\
{\rm(ii)} For any $s\in\R$ and $\ell\in\{1,\dotsc,n\}$, let 
$u\in H^s_{\ov\Omega}(\R^n,\Wedge^\ell)$ satisfy $du=0$ in $\R^n$, and $\int\!u=0$ if $\ell=n$. Then there exists $v\in H^{s+1}_{\ov\Omega}(\R^n,\Wedge^{\ell-1})$ such that $dv=u$, and there is a constant $C$ independent of $u$ such that
$$
   \Vert v\Vert_{H^{s+1}(\R^n)} \le C\, \Vert u\Vert_{H^{s}(\R^n)}\;.
$$
\end{proposition}

\begin{proof}
With $du=0$ ($\underline du=0$ in case (ii)), the relations \eqref{eq:R&T&d} reduce to 
$$
 u=d\,R_\ell u \quad\text{ and }\quad u=d\,T_\ell u\,. 
$$
Therefore in case (i) we take $v= R_\ell u$ and in case (ii) $v= T_\ell u$. The estimates are a consequence of the boundedness of the operators $R_\ell$ and $T_\ell$ as given in Corollary~\ref{C:QRSTH^s}.
\end{proof}

In the case $s=0$, there is a natural isomorphism (extension by zero outside $\Omega$) between the spaces $L^2(\Omega,\Wedge^\ell)$  and $L^2_{\ov\Omega}(\R^n,\Wedge^\ell)$. Thus for a differential form $u\in L^2(\Omega,\Wedge^\ell)$, both (i) and (ii) of the Proposition can be applied, giving a solution $v$ of $dv=u$ with coefficients in $H^1(\Omega)$ for case (i) and -- apparently stronger -- in $H^1_0(\Omega)$ for case (ii). It is important to notice, however, that the condition $du=0$ does not mean the same thing in both cases:

In case (i), it simply means $du=0$ in the sense of distributions in the open set $\Omega$. In case (ii), the condition is $du=0$ in the sense of distributions on $\R^n$, and this is stronger: It includes not only $du=0$ inside $\Omega$, but also a boundary condition $\nu\wedg u=0$ on $\partial\Omega$ in a weak sense. 

\subsection{Differential forms with polynomial coefficients}

As we have seen, the Poincar\'e-type operator $R_\ell$ preserves the class of differential forms with polynomial coefficients. This class has recently attracted some attention in the field of finite element methods. For quite a while already in relation with numerical methods for electromagnetism \cite{Hiptmair02}, but more recently also in other applications including elasticity theory \cite{ArnoldFalkWinther06}, finite dimensional subcomplexes of the de Rham complex generated by polynomials have been studied. 

For the following, we assume we have a piece of such a complex, namely for some $\ell\in\{1,\dotsc,n\}$ two spaces $P(\Wedge^{\ell-1})$ and $P(\Wedge^{\ell})$ of differential forms of order $\ell-1$ and $\ell$ with coefficients which are polynomials in $x_1,\dotsc,x_n$, which we require to satisfy the following two conditions:\\
\textbf{1. } The space $P(\Wedge^{\ell})$ is invariant with respect to dilations and translations, that is
$$
   \text{For any }t\in\R, a\in\R^n:
   \text{ If }u\in P(\Wedge^{\ell}), \text{ then }
   \bigl(x\mapsto u(tx+a)\bigr) \in P(\Wedge^{\ell}) \;.
$$ 
\textbf{2. } The interior product (``Koszul'' multiplication) $x\contr:u\mapsto x\contr u$ maps 
$P(\Wedge^{\ell})$ to $P(\Wedge^{\ell-1})$.

Then, as in Section~\ref{S:BogoPoin}, we fix a function $\theta\in\Cinf_0(\R^n)$ with support in a ball $B$ satisfying $\int\!\theta(x)\,dx=1$, and we define the Poincar\'e-type operator $R_\ell$ as in Definition~\ref{D:4iops}.

\begin{proposition}
\label{P:Rlpoly}
 The operator $R_\ell$ maps $P(\Wedge^{\ell})$ into $P(\Wedge^{\ell-1})$, and for any bounded domain $\Omega\subset\R^n$ that is starlike with respect to the ball $B$ and for any $s\in\R$ there is a constant $C$ such that for all $u\in P(\Wedge^{\ell})$
$$
  \Vert R_\ell u\Vert_{H^{s+1}(\Omega)} \,\le\, C\,\Vert u \Vert_{H^{s}(\Omega)}\;.
$$
In addition, we have for all $u\in P(\Wedge^{\ell})$
$$
   u \,=\, d\,R_\ell u \,+\, R_{\ell+1}\,du \;.
$$
\end{proposition}
\begin{proof}
 That $R_\ell$ maps $P(\Wedge^{\ell})$ into $P(\Wedge^{\ell-1})$ is a consequence of the representation \eqref{eq:regPoincare} and conditions 1. and 2. The estimate follows from the continuity stated in  Corollary~\ref{C:QRSTH^s}. 
\end{proof}

In \cite{ArnoldFalkWinther06}, complexes of polynomial differential forms are studied that satisfy conditions 1. and 2. above, and in fact a more restrictive condition than 1., namely invariance with respect to all affine transformations. The latter condition is suitable for finite elements on simplicial meshes, but our more general condition 1. covers also some cases of polynomials used in finite elements on tensor product meshes. A well-known example in 3 dimensions is the complex studied for example in \cite{CoDaDe08}, which uses spaces $Q^{p_1,p_2,p_3}$ of polynomials of partial degree $p_j$ in the variable $x_j$, $j=1,2,3$. The complex is then for a given $p\in\N$
$$
 P(\Wedge^0) \to^d
 P(\Wedge^1) \to^d
 P(\Wedge^2) \to^d
 P(\Wedge^3) 
$$
with
\begin{align*}
  P(\Wedge^0) &= Q^{p,p,p}(\Wedge^0)\;,\\
  P(\Wedge^1) &=
  \bigl\{u_1dx_1+u_2dx_2+u_3dx_3 \mid
    u_1\in Q^{p-1,p,p},\;
    u_2\in Q^{p,p-1,p},\;
    u_3\in Q^{p,p,p-1} \bigr\}\;,\\
  P(\Wedge^2) &=
    \bigl\{u_1dx_2\wedg dx_3+u_2dx_3\wedg dx_1+u_3dx_1\wedg dx_2 \mid \\
    &\hphantom{=\bigl\{u_1dx_2\wedg dx_3+u_2dx}
     u_1\in Q^{p,p-1,p-1},\;
    u_2\in Q^{p-1,p,p-1},\;
    u_3\in Q^{p-1,p-1,p} \bigr\}\;,\\
  P(\Wedge^3) &= Q^{p-1,p-1,p-1}(\Wedge^3)\;.
\end{align*}
It is clear that these spaces form a subcomplex of the de Rham complex, and that they satisfy conditions 1. and 2. above.

\subsection{Bounded Lipschitz domains}

In this subsection we draw some conclusions from Theorem~\ref{T:pseudo} that are valid for bounded Lipschitz domains. The main property of a bounded Lipschitz domain $\Omega$ that is relevant here is the existence of a finite covering of $\ov\Omega$ by open sets $U_i$, $i=1,\dots,m$ such that each $U_i\cap\Omega$ is starlike with respect to a ball $B_i$, and a subordinate partition of unity $(\chi_i)_{i=1,\dots,m}$. 
This means that 
 $\chi_i\in \Cinf_0(\R^n)$, $\supp\chi_i\subset U_i$, and
 $\sum_{i=1}^m \chi_i(x)=1$ for all $x$ in a neighborhood of $\overline\Omega$.
 
For each $i=1,\dots,m$ we can choose a smoothing function $\theta_i$ supported in $B_i$ and satisfying $\int\!\theta_i(x)dx=1$ and define the integral operators $R_{\ell,i}$ and $T_{\ell,i}$ accordingly. By Theorem~\ref{T:pseudo}, these are all pseudodifferential operators of order $-1$ on $\R^n$. They all satisfy the homotopy relations \eqref{eq:R&T&d}, but they do not have good support properties with respect to $\Omega$, only with respect to their respective $U_i\cap\Omega$.
We then define operators $R_\ell$ and $T_\ell$ according to
\begin{equation}
\label{eq:RlTlonLip}
 R_\ell u = \sum_{i=1}^m \chi_i R_{\ell,i} u\quad\text{and }
 T_\ell u = \sum_{i=1}^m T_{\ell,i}(\chi_i u)\quad\text{for }
 u\in\Cinf_0(\R^n,\Wedge^\ell), \; 1\le\ell\le n\,.
\end{equation}
These operators are still pseudodifferential operators of order $-1$ on $\R^n$, but they have better support properties with respect to $\Omega$:

If $u\in\Cinf_0(\R^n,\Wedge^\ell)$ vanishes in $\Omega$, then it vanishes in $U_i\cap\Omega$, and since $U_i\cap\Omega$ is starlike with respect to $B_i$, 
$R_{\ell,i}u$ vanishes in $U_i\cap\Omega$ and therefore $\chi_i R_{\ell,i}u$ vanishes in all of $\Omega$. Hence $R_\ell u$ vanishes in $\Omega$. In other words, the restriction of $R_\ell u$ to $\Omega$ depends only on the restriction of $u$ to $\Omega$.

For $T_\ell$ the argument is similar: If $\supp u\subset\ov\Omega$, then 
$\supp \chi_i u\subset\ov{U_i\cap\Omega}$, and therefore \\
$\supp T_{\ell,i}(\chi_i u)\subset\ov{U_i\cap\Omega}\subset\ov\Omega$. Hence 
$\supp T_\ell u\subset\ov\Omega$.

As a result, we immediately get the same mapping properties as in 
Corollaries~\ref{C:QRSTCinf} and
\ref{C:QRSTH^s}.
\begin{lemma}
\label{L:RTmappingLip}
 Let $\Omega\subset\R^n$ be a bounded Lipschitz domain and let the operators $R_\ell$ and $T_\ell$ for $1\le\ell\le n$ be defined from a finite starlike open cover of $\ov\Omega$ as in \eqref{eq:RlTlonLip}. 
Then $R_\ell$ defines continuous mappings
from $\Cinf(\Omega,\Wedge^\ell)$ to $\Cinf(\Omega,\Wedge^{\ell-1})$,
from $\Cinf(\ov\Omega,\Wedge^\ell)$ to $\Cinf(\ov\Omega,\Wedge^{\ell-1})$,
and for any $s\in\R$ 
from $H^s(\Omega,\Wedge^\ell)$ to $H^{s+1}(\Omega,\Wedge^{\ell-1})$.
The operator $T_\ell$ defines continuous mappings
from $\Cinf_0(\Omega,\Wedge^\ell)$ to $\Cinf_0(\Omega,\Wedge^{\ell-1})$,
from $\Cinf_{\ov\Omega}(\R^n,\Wedge^\ell)$ to $\Cinf_{\ov\Omega}(\R^n,\Wedge^{\ell-1})$,
and for any $s\in\R$ 
from $H^s_{\ov\Omega}(\R^n,\Wedge^\ell)$ to $H^{s+1}_{\ov\Omega}(\R^n,\Wedge^{\ell-1})$.
\end{lemma}

On the other hand, the simple anticommutation relations \eqref{eq:R&T&d} are, of course, no longer valid for these composite operators $R_\ell$ and $T_\ell$. Instead we have for $1\le\ell\le n-1$
$$
\begin{aligned}
  \bigl(dR_\ell + R_{\ell+1}d\bigr) u 
  &= 
   d\sum_{i=1}^m \chi_i R_{\ell,i} u + \sum_{i=1}^m \chi_i R_{\ell+1,i}du\\
  &=
   \sum_{i=1}^m \chi_i\bigl(dR_{\ell,i}+R_{\ell+1,i}d\bigr) u +
   \sum_{i=1}^m [d,\chi_i]R_{\ell,i} u\\
  &=
   \sum_{i=1}^m \chi_i u - K_{\ell} u
   \quad\text{ with }\;
   K_{\ell} u = -\sum_{i=1}^m [d,\chi_i]R_{\ell,i} u\;.
\end{aligned}
$$
On a neighborhood of $\overline\Omega$, this reduces to
$\quad
 \bigl(dR_\ell + R_{\ell+1}d\bigr) u = u - K_{\ell} u\;.
$

From the product rule $d(\chi_i u)= (d\chi_i)\wedg u + \chi_i du$
we obtain the commutator $[d,\chi_i]u=(d\chi_i)\wedg u$, and hence the expression for $K_{\ell}$:
\begin{equation}
\label{eq:defKl}
  K_{\ell} u = -\sum_{i=1}^m (d\chi_i)\wedg R_{\ell,i}u,\quad 1\le\ell\le n\;.
\end{equation}
This shows immediately that $K_{\ell}$ is a pseudodifferential operator of order $-1$ on $\R^n$, and that it has the same support properties as the operator $R_\ell$.

To complete the family for the endpoints $\ell=0$ and $\ell=n$, we notice that for a $0$-form $u$
$$
  R_1 du = \sum_{i=1}^m \chi_i R_{1,i}du =
    \sum_{i=1}^m \chi_i \bigl(u - (\theta_i,u)\bigr)
$$  
and for an $n$-form $u$
$$
  dR_nu = d \sum_{i=1}^m \chi_i R_{n,i} u = 
   \sum_{i=1}^m \chi_i dR_{n,i} u +
   \sum_{i=1}^m [d,\chi_i]R_{n,i} u =
   \sum_{i=1}^m \chi_i u +
   \sum_{i=1}^m d\chi_i\wedg R_{n,i} u
$$
Therefore if we set $H^s(\Omega,\Wedge^{-1}) = H^s(\Omega,\Wedge^{n+1}) = \{0\}$,
$R_0u=0$, $K_0=\sum_{i=1}^m (\theta_i,u)\chi_i$,  $R_{n+1}=0$, we obtain the homotopy relation for the de Rham complex without boundary conditions \eqref{eq:dRwobc}
\begin{equation}
\label{eq:dR+Rd}
   d\,R_\ell u + R_{\ell+1} \, d u  = u -K_\ell u\quad \text{ for all }\;
  0\le\ell\le n\,.
\end{equation}
Note that this relation is now valid only in a neighborhood of $\ov\Omega$, not in all of $\R^n$.
As a consequence of \eqref{eq:dR+Rd} we get 
$$
  d K_\ell u = d u - d R_{\ell+1}du = K_{\ell+1}d u \quad\text{ for all  $0\le\ell\le n$}\;.
$$  

For the operator $T_\ell$ we obtain similarly, when $1\leq \ell \leq n-1$,
$$
  \bigl(dT_\ell + T_{\ell+1}d\bigr) u 
  =
   (\sum_{i=1}^m \chi_i) u - L_{\ell} u
   \quad\text{ with }
   L_{\ell} u = \sum_{i=1}^m T_{\ell+1,i}[d,\chi_i] u \;.
$$
On a neighborhood of $\overline\Omega$, this reduces to
$ \bigl(dT_\ell + T_{\ell+1}d\bigr) u = u - L_{\ell} u$
with the pseudodifferential operator $L_{\ell}$ of order $-1$ given by
\begin{equation}
\label{eq:defLl}
  L_{\ell} u = \sum_{i=1}^m T_{\ell+1,i}((d\chi_i)\wedg u),\quad 0\le\ell\le n-1\;.
\end{equation}
We complete this with 
$H^{s}_{\ov\Omega}(\R^n,\Wedge^{-1}) = H^{s}_{\ov\Omega}(\R^n,\Wedge^{n+1})=\{0\}$,
$T_0=0$, $T_{n+1}=0$, and $L_nu=\sum (\int\chi_i u)\star\theta_i$
and obtain the homotopy relation for the de Rham complex with compact support \eqref{eq:dRwcs}
\begin{equation}
\label{eq:dT+Td}
  d\,T_\ell u + T_{\ell+1} \, du = u-L_\ell u \quad \text{ for all }\;
  0\le\ell\le n.
\end{equation}
This relation is valid in a neighborhood of $\ov\Omega$, but now if we apply it to a $u$ with support in $\ov\Omega$, it will be valid in all of $\R^n$. Again as before we obtain 
$$
   d L_\ell u = L_{\ell+1} d u \quad\text{ for all  $0\le\ell\le n$}\;. 
$$
\begin{Remark}
In this subsection on Lipschitz domains, we are using the extended de Rham complexes \eqref{eq:dRwobc} and \eqref{eq:dRwcs}, rather than the sequences \eqref{eq:edRwobc} and \eqref{eq:edRwcs} as we did for starlike domains. For this reason, we now have $R_0=0$, $T_0=0$, $R_{n+1}=0$ and $T_{n+1}=0$.
\end{Remark}

Before drawing conclusions, we prove a stronger version of the relations 
\eqref{eq:dR+Rd} and \eqref{eq:dT+Td}, where the perturbations of the identity $K_\ell$ and $L_\ell$ are not just of order $-1$, but in fact infinitely smoothing in a neighborhood of $\ov\Omega$. 

Let $x_0\in\R^n$. We shall say that the family of functions $(\chi_i)_{i=1,\dots,m}$ is \emph{flat at $x_0$} if each $\chi_i$ is constant in a neighborhood of $x_0$. 
We will also call an open covering $(U_i)_{i=1,\dots,m}$ of $\ov\Omega$ by a slight abuse of language \emph{starlike} if each $U_i\cap\Omega$ is starlike with respect to some open ball $B_i$.
\begin{lemma}
\label{L:flatPoU}
Let $\Omega$ be a bounded Lipschitz domain.
Then there exists a finite number of starlike finite open coverings $(U_i^{(j)})_{i=1,\dots,m^{(j)}}$, $j=1,\dots,k$, of $\ov\Omega$ and subordinate partitions of unity, such that for any $x_0\in\R^n$ at least one of the partitions of unity is flat at $x_0$.
\end{lemma}
\begin{proof}
In a first step we show that for a given $x_0\in\R^n$ there exists a starlike finite open covering $(U_i)_{i=0,\dots,m}$ of $\overline\Omega$  and a partition of unity subordinate to this covering which is flat at $x_0$.

Let first $x_0\in\overline\Omega$. Let $U_0$ be a neighborhood of $x_0$ such that $U_0\cap\Omega$ is starlike with respect to a ball, $V_0$ another neighborhood of $x_0$ such that $\ov V_0\subset U_0$ and $\chi_0\in \Cinf_0(\R^n)$ such that $\supp\chi_0\subset U_0$ and $\chi_0\equiv1$ on a neighborhood of $\ov V_0$. We may assume that $\Omega\setminus\ov V_0$ is still Lipschitz. 
Choose a finite open covering $(U_i)_{i=1,\dots,m}$ of $\overline{\Omega\setminus V_0}$ such that each $U_i\cap\Omega$ is starlike with respect to a ball. 
Let $\{\tilde\chi_i\mid i=1,\cdots,m\}$ be a subordinate partition of unity which therefore satisfies 
$$
  \sum_{i=1}^m \tilde\chi_i(x)=1 \quad
  \text{ for all $x$ in a neighborhood of } \overline{\Omega\setminus V_0}\,.
$$
Then defining for $i=1,\cdots,m$:
$$
  \chi_i = (1-\chi_0)\tilde\chi_i 
$$
we have a starlike covering $(U_i)_{i=0,\dots,m}$ of $\overline\Omega$ and a subordinate partition of unity $(\chi_i)_{i=0,\dots,m}$ which is flat at $x_0$.

If now $x_0\in\R^n\setminus\ov\Omega$, then from any partition of unity subordinate to an open covering of $\ov\Omega$ we get another one which is flat at $x_0$ by multiplying with a cut-off function which is $1$ on a neighborhood of $\ov\Omega$ and vanishes on a neighborhood of $x_0$. 

In a second step we choose $R>0$ such that $\ov\Omega\subset B_R(0)$.
To any $x_0\in\ov B_R(0)$ there exists, as we have proved in the first step, a neighborhood $V(x_0)$ and a starlike open covering $(U_i^{(x_0)})_i$ of $\ov\Omega$ with a subordinate partition of unity $(\chi_i^{(x_0)})_i$ which is flat at any point of $V(x_0)$. 
The open covering $\bigl(V(x_0)\bigr)_{x_0\in\overline B_R(0)}$ of the compact set $\overline B_R(0)$ contains a finite subcovering associated with points 
$x_0=x_1,\dots,x_k\in\overline B_R(0)$. The corresponding family of open coverings 
$(U_i^{(x_j)})$ and partitions of unity $(\chi_i^{(x_j)})$ for $j=1,\dots,k$ will have the required properties for all points $x_0\in\overline B_R(0)$. For the remaining points $x_0\in\R^n\setminus\overline B_R(0)$, one adds one of the previous partitions of unity, after multiplying each of its functions by a $\Cinf$ cut-off function that is $1$ in a neighborhood of $\overline\Omega$ and has its support in $B_R(0)$.
\end{proof}

\begin{theorem}
\label{T:homCinf}
Let $\Omega\subset\R^n$ be a bounded Lipschitz domain. Then for $\ell=0,1,\dotsc,n$, there exist pseudodifferential operators $R_\ell$, $T_\ell$ of order $-1$ and $K_\ell$, $L_\ell$ of order $-\infty$ on $\R^n$ with the following properties:\\
{\rm(i)}
The operators define continuous mappings
\begin{align*}
  R_\ell&: \Cinf(\ov\Omega,\Wedge^\ell) \to \Cinf(\ov\Omega,\Wedge^{\ell-1})
  &\text{ and }\quad
  T_\ell&: \Cinf_{\ov\Omega}(\R^n,\Wedge^\ell) 
       \to \Cinf_{\ov\Omega}(\R^n,\Wedge^{\ell-1})\,,\\
\intertext{
and for any $s\in\R$}
  R_\ell&: H^s(\Omega,\Wedge^\ell) \to H^{s+1}(\Omega,\Wedge^{\ell-1})
  &\text{ and }\quad
  T_\ell&: H^s_{\ov\Omega}(\R^n,\Wedge^\ell) 
       \to H^{s+1}_{\ov\Omega}(\R^n,\Wedge^{\ell-1})\,,\\
  K_\ell&: H^s(\Omega,\Wedge^\ell) \to \Cinf(\ov\Omega,\Wedge^{\ell})
  &\text{ and }\quad
  L_\ell&: H^s_{\ov\Omega}(\R^n,\Wedge^\ell) \to \Cinf_{\ov\Omega}(\R^n,\Wedge^{\ell})\;.
\end{align*}
{\rm(ii)} 
On a neighborhood of $\ov\Omega$, there holds for $\ell=0,1,\dotsc,n$ and any $\ell$-form $u$ on $\R^n$ with compact support
\begin{equation}
\label{dRTdKL}
  d\,R_\ell u + R_{\ell+1} \, du = u -K_\ell u
  \quad\text{ and }\quad
  d\,T_\ell u + T_{\ell+1} \, du = u -L_\ell u \;.
\end{equation}
{\rm(iii)} In particular,
$K_0$ is a finite-dimensional operator mapping $H^s(\Omega,\Lambda^0)$ continuously to $\Cinf(\ov\Omega,\Lambda^0)$ for any $s\in\R$, $L_n$ is a finite-dimensional operator mapping 
$H^s_{\ov\Omega}(\R^n,\Lambda^n)$ continuously to $\Cinf_{\ov\Omega}(\R^n,\Lambda^n)$ for any $s\in\R$, and one has in a neighborhood of $\ov\Omega$\;:
\begin{align*}
  R_1\,du &= u-K_0u\,,\quad
  T_1\,du = u-L_0u\,,\quad \text{when } \ell=0\,,\\ 
  d\,R_n u &= u - K_n u\,,\quad
  d\,T_n u = u - L_n u\,,\quad \text{when } \ell=n\,.
  \end{align*}
\end{theorem}

\begin{proof}
We give the details of the proof for the Poincar\'e-type operators $R_\ell$. For the Bogovski\u{\i}-type operators $T_\ell$, the proof is the same.

The crucial observation is that in the definitions \eqref{eq:defKl} of the perturbation operator $K_\ell$ and \eqref{eq:defLl} of $L_\ell$, the factors $d\chi_i$ are all zero in a neighborhood of any point $x_0$ in which the partition of unity $(\chi_i)_{i=1,\dotsc,m}$ is flat. The images $K_\ell u$ and $L_\ell u$ are therefore $\Cinf$ in the neighborhood of such a point (in fact, $K_\ell u$ is even zero there).

We choose now a finite number of starlike finite open coverings
$(U_i^{(j)})_{i=1,\dots,m^{(j)}}$,$j=1,\dots,k$, of $\ov\Omega$ and subordinate partitions of unity $(\chi_i^{(j)})_{i=1,\dots,m^{(j)}}$, $j=1,\dots,k$ which exist according to Lemma~\ref{L:flatPoU} in such a way that for any $x_0\in\R^n$ at least one of the partitions of unity is flat at $x_0$. For each $j=1,\dots,k$, we construct the operators $R_\ell^{(j)}$ and $K_\ell^{(j)}$ associated with the corresponding partition of unity. They satisfy the equivalent of \eqref{eq:dR+Rd} on a neighborhood of $\ov\Omega$, namely
\begin{equation}
\begin{aligned}
\label{eq:RKcomm}
  \bigl(d\,R_\ell^{(j)} + R_{\ell+1}^{(j)}\,d\,\bigr) u &= u - K_{\ell}^{(j)} u\,,\\
  dK_{\ell}^{(j)} u &= K_{\ell+1}^{(j)} \,du\,.
\end{aligned}
\end{equation}
We can then define
$$
\begin{aligned}
   R_\ell &= R_\ell^{(1)} + K_{\ell-1}^{(1)}R_\ell^{(2)} + 
                K_{\ell-1}^{(1)}K_{\ell-1}^{(2)}R_\ell^{(3)} +\dots+ 
                K_{\ell-1}^{(1)}\cdots K_{\ell-1}^{(k-1)} R_\ell^{(k)}\\
   K_\ell &= K_\ell^{(1)}\cdots K_\ell^{(k)} \;.
\end{aligned}
$$
Using the relations \eqref{eq:RKcomm}, one can easily verify 
that on a neighborhood of $\overline\Omega$ we have
\begin{equation}
\label{eq:hatRK}
 \bigl(d\, R_\ell +  R_{\ell+1}\,d\,\bigr) u = u - K_{\ell} u
 \quad\text{ and }\quad d K_\ell u = K_{\ell+1} d u\;.
\end{equation}

In addition, we find that the operator $K_{\ell}$ is not only a pseudodifferential operator of order $-k$ as a product of pseudodifferential operators of order $-1$, but actually of order $-\infty$, that is, an integral operator with $\Cinf$ kernel, continuously mapping ${\cal D}'(\R^n)$ to $\Cinf(\R^n)$. The reason for this is that for any $x_0\in\R^n$, at least one of the partitions of unity $(\chi_i^{(j)})_{i=1,\dots,m^{(j)}}$ is flat at $x_0$, and that therefore the corresponding factor $K_\ell^{(j)}$ maps to functions which are $\Cinf$ in a neighborhood of $x_0$. The other factors in the definition of $K_\ell$ are pseudodifferential operators, hence pseudo-local, and therefore the product $K_\ell$ maps to functions that are $\Cinf$ in a neighborhood of $x_0$, too.
\end{proof}

The relations \eqref{dRTdKL} imply regularity results for the $d$ operator. These can be expressed as existence of solutions of maximal regularity if the solvability conditions are satisfied. We consider this first for the inhomogeneous equation $dv=u$ and then for the homogeneous equation $du=0$. Finally we obtain a regularity result for the cohomology spaces of the two de Rham complexes \eqref{eq:dRwobc} and \eqref{eq:dRwcs}.

\begin{corollary}
\label{c:Regford}
 Let $\Omega$ be a bounded Lipschitz domain in $\R^n$.\\
 For $1\le\ell\le n$ and
 any $s,t\in\R$ we have:\\
{\rm(a) } 
If $u\in H^{s}(\Omega,\Wedge^{\ell})$ satisfies $u=dv$ for some 
$v\in H^{t}(\Omega,\Wedge^{\ell-1})$,
then there exists $w\in H^{s+1}(\Omega,\Wedge^{\ell-1})$ such that $u=dw$,
and there is a constant $C$ independent of $u$ and $v$ with
$$
  \Vert w\Vert_{H^{s+1}(\Omega)} 
  \le\,C\, \bigl(
  \Vert u\Vert_{H^{s}(\Omega)} 
  + \Vert v \Vert_{H^{t}(\Omega)}
  \bigr)\;.$$
{\rm(b) } If
$u\in H^{s}_{\ov\Omega}(\R^n,\Wedge^{\ell})$ satisfies $u=dv$ for some  
$v\in H^{t}_{\ov\Omega}(\R^n,\Wedge^{\ell-1})$, 
then there exists
$w\in H^{s+1}_{\ov\Omega}(\R^n,\Wedge^{\ell-1})$ such that $u=dw$,
and there is a constant $C$ independent of $u$ and $v$ with
$$
  \Vert w\Vert_{H^{s+1}(\R^n)} 
  \le \,C\, \bigl(
  \Vert u\Vert_{H^{s}(\R^n)} 
  + \Vert v \Vert_{H^{t}(\R^n)}
  \bigr)\;.
$$
\end{corollary}
\begin{proof}
(a) If $u=dv$, then with $v=dR_{\ell-1}v+R_\ell dv +K_{\ell-1}v$ we get
$u=d\,\bigl(R_\ell u +K_{\ell-1}v\bigr)$, and $w=R_\ell u +K_{\ell-1}v$ belongs to 
$H^{s+1}(\Omega,\Wedge^{\ell-1})$ if $u\in H^{s}(\Omega,\Wedge^{\ell})$.
The estimate follows from the fact that $R_\ell$ is of order $-1$ and that $K_{\ell-1}$ maps $H^{t}(\Omega,\Wedge^{\ell-1})$ continuously to $H^{s+1}(\Omega,\Wedge^{\ell-1})$ for any $s$ and $t$.

(b) Likewise, $u=dv$ implies $u=dw$ with 
$w=T_\ell u +L_{\ell-1}v\in H^{s+1}_{\ov\Omega}(\R^n,\Wedge^{\ell-1})$ if
$u\in H^{s}_{\ov\Omega}(\R^n,\Wedge^{\ell})$.
\end{proof}

Next we consider the special case of relations \eqref{dRTdKL} where $du=0$. \begin{corollary}
\label{c:InvdCinf}
 Let $\Omega$ be a bounded Lipschitz domain in $\R^n$.
 For any $s\in\R$ and $1\le\ell\le n$ we have:\\[1ex]
{\rm(a)} 
$\quad
  u\in H^s(\Omega,\Wedge^\ell), \;du=0 \;\text{ in }\,\Omega
  \quad\Longrightarrow\quad
   u = dR_\ell u + K_{\ell}u
 \quad\text{ in }\;\Omega
$

Here $R_\ell u\in H^{s+1}(\Omega,\Wedge^{\ell-1})$
and $K_\ell u\in \Cinf(\ov\Omega,\Wedge^{\ell})$. \\[1ex]
{\rm(b)}
$\quad
  u\in H^s_{\ov\Omega}(\R^n,\Wedge^\ell), \; du=0 \;\text{ in }\,\R^n
  \quad\Longrightarrow\quad
   u = dT_\ell u + L_{\ell}u
 \quad\text{ in }\;\R^n
$

Here $T_\ell u\in H^{s+1}_{\ov\Omega}(\R^n,\Wedge^{\ell-1})$
and $L_\ell u\in \Cinf_{\ov\Omega}(\R^n,\Wedge^{\ell})$.
\end{corollary}

For a bounded Lipschitz domain $\Omega$, we consider now the cohomology spaces of regularity $s$ of the two de Rham complexes, without boundary conditions \eqref{eq:dRwobc}, and with compact support \eqref{eq:dRwcs}. Thus we introduce the corresponding two variants of the cohomology spaces, \emph{without boundary conditions}
\begin{equation}
 {\cal H}^s_\ell(\Omega) :=
  \frac{\Nullspace \bigl(d:H^{s}(\Omega,\Wedge^{\ell})\to H^{s-1}(\Omega,\Wedge^{\ell+1})\bigr)}%
  {\Range \bigl(d:H^{s+1}(\Omega,\Wedge^{\ell-1})\to H^{s}(\Omega,\Wedge^{\ell})\bigr)}
\end{equation}
and 
\emph{with compact support}:
\begin{equation}
\label{eq:Cwcs}
 {\cal H}^s_{\ov\Omega,\ell}(\R^n) 
   = \frac{\Nullspace \bigl(d:H^{s}_{\ov\Omega}(\R^n,\Wedge^{\ell})\to H^{s-1}_{\ov\Omega}(\R^n,\Wedge^{\ell+1})\bigr)}%
    {\Range \bigl(d:H^{s+1}_{\ov\Omega}(\R^n,\Wedge^{\ell-1})\to H^{s}_{\ov\Omega}(\R^n,\Wedge^{\ell})\bigr)} \;.
\end{equation}
Here we can consider the full range $0\le\ell\le n$, if we complete the complexes by $0$ as we did in \eqref{eq:dRwobc} and \eqref{eq:dRwcs}.

\begin{theorem}  
\label{T:CohomCinf}
 Let $\Omega$ be a bounded Lipschitz domain in $\R^n$, and let $0\le\ell\le n$.
 \\
{\rm(a)} For any $s\in\R$, the exterior derivatives
$$
  d:H^{s+1}(\Omega,\Wedge^{\ell-1})\to H^{s}(\Omega,\Wedge^{\ell})
  \qquad\text{ and }\quad
  d:H^{s+1}_{\ov\Omega}(\R^n,\Wedge^{\ell-1})\to H^{s}_{\ov\Omega}(\R^n,\Wedge^{\ell})
$$
define bounded operators with closed range $dH^{s+1}(\Omega,\Wedge^{\ell-1})$ and $dH^{s+1}_{\ov\Omega}(\R^n,\Wedge^{\ell-1})$\;.\\
{\rm(b)} The dimension of ${\cal H}^s_\ell(\Omega)$ is a finite number $b_\ell$ independent of $s\in\R$. Moreover
there is a $b_\ell$-dimensional subspace ${\cal H}_\ell(\ov\Omega)$ of 
$\Cinf(\ov\Omega,\Wedge^\ell)$ such that, for all $s\in\R$,
\begin{equation}\label{eq:direct}
 \Nullspace \Bigl(d:H^{s}(\Omega,\Wedge^{\ell})\to H^{s-1}(\Omega,\Wedge^{\ell+1})\Bigr)
 \, = \,
 d\, H^{s+1}(\Omega,\Wedge^{\ell-1}) \, \oplus \, {\cal H}_\ell(\ov\Omega)
  \:. 
\end{equation}
That is, for any 
$u\in H^{s}(\Omega,\Wedge^{\ell})$ satisfying $du=0$ in $\Omega$, there exist
$v\in H^{s+1}(\Omega,\Wedge^{\ell-1})$ and a unique
$w\in {\cal H}_\ell(\ov\Omega)$,  such that
$$
  u = dv + w\quad\text{with}\quad \|v\|_{H^{s+1}(\Omega)}+\|w\|_{H^s(\Omega)}\leq C_s\|u\|_{H^s(\Omega)}\;.
$$
{\rm(c)} The dimension of ${\cal H}^s_{\ov\Omega,\ell}(\R^n)$ is a finite number $\tilde b_\ell$ independent of $s\in\R$. Moreover
there is a $\tilde b_\ell$-dimensional subspace ${\cal H}_{\ov\Omega,\ell}(\R^n)$ of 
$\Cinf_{\ov\Omega}(\R^n,\Wedge^\ell)$ such that, for all $s\in\R$,
\begin{equation}\label{eq:direct1}
 \Nullspace \Bigl(d:H^{s}_{\ov\Omega}(\R^n,\Wedge^{\ell})\to H^{s-1}_{\ov\Omega}(\R^n,\Wedge^{\ell+1})\Bigr) 
 \, = \,
 d\, H^{s+1}_{\ov\Omega}(\R^n,\Wedge^{\ell-1}) \, \oplus \, {\cal H}_{\ov\Omega,\ell}(\R^n)
  \:. 
 \end{equation}
That is, for any 
$u\in H^{s}_{\ov\Omega}(\R^n,\Wedge^{\ell})$ satisfying $du=0$ in $\R^n$, there exists
$v\in H^{s+1}_{\ov\Omega}(\R^n,\Wedge^{\ell-1})$ and a unique
$w\in {\cal H}_{\ov\Omega,\ell}(\R^n)$,  such that
$$
  u = dv + w\quad\text{with}\quad \|v\|_{H^{s+1}(\R^n)}+\|w\|_{H^s(\R^n)}\leq C_s\|u\|_{H^s(\R^n)}\;.
$$
 {\rm(d)} The dimensions $b_\ell$ and $\tilde b_\ell$ are related by
$$\tilde b_{n-\ell}=b_\ell \:.$$
\end{theorem}
\begin{proof}
We give the proof for the case without boundary conditions. The proof for the case with compact support is similar if one takes into account the mapping properties of the operators $T_\ell$ and $L_\ell$.

Fix $\ell\in\{0,\dotsc,n\}$. For $s\in\R$, define 
$$
  N^s_\ell = \Nullspace \bigl(d:H^{s}(\Omega,\Wedge^{\ell})\to H^{s-1}(\Omega,\Wedge^{\ell+1})\bigr)\;
$$ 
with in particular, $N^s_n = H^{s}(\Omega, \Wedge^{n})$.
This is a closed subspace of 
$H^{s}(\Omega,\Wedge^{\ell})$, and for the study of the range of $d$, we can replace 
$H^{s+1}(\Omega,\Wedge^{\ell-1})$ by the quotient space
$$
  X^{s+1}_{\ell-1} := H^{s+1}(\Omega,\Wedge^{\ell-1}) / N^{s+1}_{\ell-1}
$$
with its natural quotient norm. We will now study the properties of $d$ as a mapping
\begin{equation}
\label{eq:dXtoN}
  d : X^{s+1}_{\ell-1} \to N^s_\ell \;.
\end{equation}
We know from \eqref{eq:hatRK} that the nullspace of $d$ is an invariant subspace of the operator $K_\ell$, and $K_\ell$ is a compact operator in $N^s_\ell$. By the same token, $K_{\ell-1}$ is defined in a natural way on the quotient space $X^{s+1}_{\ell-1}$, and it is a compact operator there.

Also from \eqref{eq:hatRK} follows that for 
$u\in H^{s+1}(\Omega,\Wedge^{\ell-1})$ we have
$$
  R_\ell du = u - K_{\ell-1}u - d R_{\ell-1}u
  \quad \equiv u - K_{\ell-1}u \mod N^{s+1}_{\ell-1} \;,
$$
and for $v\in N^s_\ell$ we have
$$
  d R_\ell v= v - K_{\ell} v\;.
$$
Together, this means that if we consider $R_\ell$ as a bounded operator from
$N^s_\ell$ to $X^{s+1}_{\ell-1}$, it defines a two-sided regularizer (inverse modulo compact operators) of the operator $d$ in \eqref{eq:dXtoN}. By the well-known theory of Fredholm operators, this implies that $d$ in \eqref{eq:dXtoN} is a Fredholm operator. Its image is therefore closed, which proves point (a), and it has finite codimension, which shows that ${\cal H}^s_\ell(\Omega)$ is finite dimensional.

Let us now define the direct summand ${\cal H}_\ell(\ov\Omega)$. Let $b_\ell = \dim\,{\cal H}^0_\ell$.
It is a consequence of the above results that $dH^1(\Omega,\Lambda^{\ell-1})$ has a $b_\ell$-dimensional direct summand, say $\tilde{\cal H}(\Omega,\Lambda^{\ell})$ in $N^0_\ell$. That is 
$$
  N^0_\ell = dH^1(\Omega,\Lambda^{\ell-1}) \oplus \tilde{\cal H}(\Omega,\Lambda^{\ell})\;.
$$ 
Define
$$
  {\cal H}_\ell(\ov\Omega) = K_\ell  \tilde{\cal H}(\Omega,\Lambda^{\ell})\subset \Cinf(\ov\Omega,\Wedge^\ell)\;.
$$
Then, by \eqref{eq:hatRK}, ${\cal H}_\ell(\ov\Omega)\subset N^s_\ell$ for all $s$. 
Moreover $dH^{s+1}(\Omega,\Lambda^{\ell-1}) \cap {\cal H}_\ell(\ov\Omega)=\{0\}$. To see this, suppose that $dv=K_\ell w$ where $v\in H^{s+1}(\Omega,\Lambda^{\ell-1})$ and $w\in \tilde{\cal H}(\Omega,\Lambda^{\ell})$. Thus, using \eqref{dRTdKL},
$d(R_\ell dv + K_{\ell-1} v)=w-dR_\ell w$ and hence $du=w$ where $u=R_\ell K_\ell w +K_{\ell-1}v +R_\ell w \in \Cinf(\ov\Omega,\Lambda^{\ell-1})\subset H^1(\Omega, \Lambda^{\ell-1})$. So, by the definition of   $\tilde{\cal H}(\Omega,\Lambda^{\ell})$, $w=0$ and then again, $dv=K_\ell w=0$. In a similar way, we can show that $K_\ell$ is one--one on ${\cal H}_\ell(\ov\Omega)$, so that $\dim {\cal H}_\ell(\ov\Omega) = b_\ell$.

We next prove \eqref{eq:direct}.
Given $u\in N^s_\ell$, write $u=dR_\ell u + K_\ell (dR_\ell u+K_\ell u)$. Now
$K_\ell u \in \Cinf(\ov\Omega,\Wedge^\ell) \subset H^0(\Omega,\Wedge^{\ell})$, so by the definition of $\tilde{\cal H}(\Omega,\Lambda^{\ell})$, we can write 
$$
  K_\ell u = dv' + w'
  \quad\text{ with }\quad
  v'\in H^{1}(\Omega,\Wedge^{\ell-1}),\,
  w'\in \tilde{\cal H}(\Omega,\Lambda^\ell) \;.
$$
Hence $u=dv+w$ with $v=R_\ell u +K_{\ell-1}R_\ell u + K_{\ell-1} v'
\in H^{s+1}(\Omega,\Wedge^{\ell-1})$, $w=K_\ell w'\in {\cal H}_\ell(\ov\Omega)$, and
$\|v\|_{H^{s+1}}+\|w\|_{H^s}\leq C_s\|u\|_{H^s}$.

It is a consequence of \eqref{eq:direct} that ${\cal H}^s_\ell(\Omega)$ is isomorphic to 
${\cal H}_\ell(\ov\Omega)$, and hence $\dim {\cal H}^s_\ell(\Omega)=b_\ell$ for all $s$.

To prove part (d), observe that
\begin{align*}
 \Bigl\{\Nullspace \Bigl(d:H^{s}(\Omega,\Wedge^{\ell})\to H^{s-1}(\Omega,\Wedge^{\ell+1})\Bigr)\Bigr\}^\perp 
&=\,\delta H^{-s+1}_{\ov\Omega}(\R^n,\Wedge^{\ell+1})\\
&=\,*dH^{-s+1}_{\ov\Omega}(\R^n,\Wedge^{n-\ell-1})
\end{align*}
and
\begin{equation*}
\bigl\{dH^{s+1}(\Omega, \Wedge^{\ell-1})\bigr\}^\perp\,=\,*\Nullspace \Bigl(d:H^{-s}_{\ov\Omega}(\R^n,\Wedge^{n-\ell})\to H^{-s-1}_{\ov\Omega}(\R^n,\Wedge^{n-\ell+1})\Bigr)\;.
\end{equation*}
Therefore, by duality, 
\begin{align*}
b_\ell &= \dim\left\{ \frac{\Nullspace \bigl(d:H^{s}(\Omega,\Wedge^{\ell})\to H^{s-1}(\Omega,\Wedge^{\ell+1})\Bigr)}
  {\Range \bigl(d:H^{s+1}(\Omega,\Wedge^{\ell-1})\to H^{s}(\Omega,\Wedge^{\ell})\bigr)}\right\}\\  &= \dim\left\{\frac{\Nullspace \bigl(d:H^{-s}_{\ov\Omega}(\R^n,\Wedge^{n-\ell})\to H^{-s-1}_{\ov\Omega}(\R^n,\Wedge^{n-\ell+1})\bigr)}
    {\Range \bigl(d:H^{-s+1}_{\ov\Omega}(\R^n,\Wedge^{n-\ell-1})\to H^{-s}_{\ov\Omega}(\R^n,\Wedge^{n-\ell})\bigr)}\right\} = \tilde b_{n-\ell}\;.
\end{align*}
\end{proof}

\begin{Remark}
When $\ell = 0$, then 
\begin{align*} {\cal H}_0(\ov\Omega)= \Nullspace \bigl(d:H^{s}(\Omega,\Wedge^{0})\to H^{s-1}(\Omega,\Wedge^{1})\bigr) &= \R \quad \text{(the constant functions) and}\\
{\cal H}_{\ov\Omega,0}(\R^n)=\Nullspace \bigl(d:H^{s}_{\ov\Omega}(\R^n,\Wedge^{0})\to H^{s-1}_{\ov\Omega}(\R^n,\Wedge^{1})\bigr)&=\{0\}\end{align*}
so, by duality, 
\begin{align*}
dH^{s+1}(\Omega,\Wedge^{n-1}) &=H^s(\Omega,\Wedge^n)\ ,\quad {\cal H}_n(\ov\Omega)=\{0\}\ , \\
dH^{s+1}_{\ov\Omega}(\R^n,\Wedge^{n-1}) &=\{u\in H^s_{\ov\Omega}(\R^n,\Wedge^n):\smallint u = 0\}\ ,
\end{align*}
and ${\cal H}_{\ov\Omega,n}(\R^n)$ can be taken to be $\{cL_n 1_{\ov\Omega}\mid c\in\R\}$ where $1_{\ov\Omega}$ is the characteristic function of $\ov\Omega$. 
Therefore
$b_0=\tilde b_n=1$ and $b_n =\tilde b_0=0$. 

When $1\leq \ell \leq n-1$, we can take $\tilde{\cal H}(\Omega,\Lambda^\ell)$ to be the orthogonal complement of $dH^1(\Omega,\Lambda^{\ell-1})$ in $N^0_\ell$, so that 
$$
 {\cal H}_\ell(\ov\Omega)= K_\ell \{u\in L^2(\Omega,\Lambda^\ell) \mid
 du=0, \delta u=0 \text { and } \nu\contr u=0 \text{ on } \partial\Omega\}
 \,, 
$$
Similarly we can take
$$
 {\cal H}_{\ov\Omega,\ell}(\R^n)= L_\ell {\cal E_\ell}\{u\in L^2(\Omega,\Lambda^\ell) \mid
 du=0, \delta u=0 \text { and } \nu\wedg u=0 \text{ on } \partial\Omega\}
$$ 
where ${\cal E}_\ell:L^2(\Omega,\Lambda^\ell)\to L^2(\R,\Lambda^\ell)$ denotes extension by zero.
 
 The integers $b_\ell$ are the Betti numbers of $\Omega$.
\end{Remark} 

Note that the sequence of Betti numbers $b_0,\dotsc,b_n$ will in general be different from the sequence $\tilde b_0,\dotsc,\tilde b_n$. For example, for the standard torus embedded in $\R^3$, one finds without difficulties the two sequences $1,1,0,0$ and $0,0,1,1$, and for the ball with a hole $B_2(0)\setminus \ov{B_1(0)}$, one gets the two sequences $1,0,1,0$ and $0,1,0,1$.

Classically, one considers the de Rham complexes for differential forms with smooth coefficients
\begin{equation}
\label{eq:CinfdRwobc}
 0 \to\Cinf(\ov\Omega,\Wedge^0) \to^d \Cinf(\ov\Omega,\Wedge^1)
 \to^d\cdots \to^d \Cinf(\ov\Omega,\Wedge^n) \to 0
\end{equation}
and 
\begin{equation}
\label{eq:CinfdRwcs}
 0 \to \Cinf_{\ov\Omega}(\R^n,\Wedge^0) \to^d \Cinf_{\ov\Omega}(\R^n,\Wedge^1) \to^d\cdots
 \to^d \Cinf_{\ov\Omega}(\R^n,\Wedge^n)  \to 0
\end{equation}
With the same arguments as in the preceding proof one can see that the associated cohomology spaces are isomorphic to those with finite regularity considered in Theorem~\ref{T:CohomCinf}. It suffices to notice that pseudodifferential operators map $\Cinf$ functions to $\Cinf$ functions.

\begin{corollary}
\label{C:CinfCohom}
 Let $\Omega$ be a bounded Lipschitz domain in $\R^n$ and $0\le\ell\le n$ .\\
{\rm(a)} The cohomology space without boundary condition
$$
 \frac{\Nullspace \bigl(d:\Cinf(\ov\Omega,\Wedge^{\ell})\to \Cinf(\ov\Omega,\Wedge^{\ell+1})\bigr)}
  {\Range \bigl(d:\Cinf(\ov\Omega,\Wedge^{\ell-1})\to \Cinf(\ov\Omega,\Wedge^{\ell})\bigr)}
$$
of the de Rham complex \eqref{eq:CinfdRwobc} has dimension $b_\ell$  and is isomorphic to 
${\cal H}_\ell(\ov\Omega)$. There is a splitting
$$
  \Nullspace \Bigl(d:\Cinf(\ov\Omega,\Wedge^{\ell})\to \Cinf(\ov\Omega,\Wedge^{\ell+1})\Bigr)
  \,=\,
  d\, \Cinf(\ov\Omega,\Wedge^{\ell-1}) \,\oplus\, {\cal H}_\ell(\ov\Omega)
  \:.
$$
{\rm(b)} The cohomology space with compact support
$$
   \frac{\Nullspace \bigl(d:\Cinf_{\ov\Omega}(\R^n,\Wedge^{\ell})\to \Cinf_{\ov\Omega}(\R^n,\Wedge^{\ell+1})\bigr)}
  {\Range \bigl(d:\Cinf_{\ov\Omega}(\R^n,\Wedge^{\ell-1})\to \Cinf_{\ov\Omega}(\R^n,\Wedge^{\ell})\bigr)}
$$
of the de Rham complex \eqref{eq:CinfdRwcs} has dimension $\tilde b_\ell$ and is isomorphic to 
${\cal H}_{\ov\Omega,\ell}(\R^n)$. There is a splitting
$$
\Nullspace \Bigl(d:\Cinf_{\ov\Omega}(\R^n,\Wedge^{\ell})\to \Cinf_{\ov\Omega}(\R^n,\Wedge^{\ell+1})\Bigr)
  \,=\,
  d\, \Cinf_{\ov\Omega}(\R^n,\Wedge^{\ell-1}) \,\oplus\, {\cal H}_{\ov\Omega,\ell}(\R^n)
  \:.
$$

\end{corollary}

\begin{Remark}
All the results of this section remain valid when $H^s$ is replaced by $B^s_{pq}\, (0<p\leq \infty, 0<q\leq \infty)$, or by $F^s_{pq}\, (0<p<\infty, 0<q\leq \infty)$. 

We make the following additional comments.

In Corollary~\ref{c:Regford}, all that is required of $v$ is that, in part (a), $v$ be the restriction to $\Omega$ of a distribution (with compact support) on $\R^n$, while in part (b), $v$ be a distribution on $\R^n$ with support in $\overline\Omega$. Indeed, it is well known that distributions with compact support are of finite order, so there exists then a finite index $t$ such that $v$ belongs to one of the spaces required in the corollary.

The dimension of the cohomology spaces ${\cal H}^s_\ell(\Omega)$ and 
${\cal H}^s_{\overline\Omega,\ell}$, defined using $B^s_{pq}$ or $F^s_{pq}$ in place of $H^s$, are still equal to $b_l$ and $\tilde b_l$.
\end{Remark}

We conclude by mentioning that we have now proved Theorem \ref{intro}, stated in the Introduction.

\medskip

{\small \textbf{Acknowledgments: }
Part of this research was conducted while the first author was invited to the Centre for Mathematics and its Applications at the Australian National University, Canberra, and to the Institute for Computational Engineering and Sciences at the University of Texas at Austin, USA.
The second author was supported by the Australian Government through the Australian Research Council. 
}

\textbf{Addresses:}\\[1ex]
\textsc{Martin Costabel}\\
IRMAR, Universit\'e de Rennes 1,
Campus de Beaulieu\\35042 Rennes Cedex, France\\
E-mail: \texttt{martin.costabel@univ-rennes1.fr}
\\[1ex]
\textsc{Alan McIntosh}\\
Centre for Mathematics and its Applications, Australian National University\\
Canberra, ACT 0200, Australia\\
E-mail: \texttt{alan@maths.anu.edu.au}

\end{document}